\documentclass[12pt,a4paper,reqno]{amsart}
\usepackage{amssymb,amsmath,latexsym,oldgerm,mathrsfs,epsfig}
\numberwithin{equation}{section}
\newtheorem{theorem}{Theorem}[section]
\newtheorem{lemma}[theorem]{Lemma}
\newtheorem{proposition}[theorem]{Proposition}
\newtheorem{corollary}[theorem]{Corollary}
\theoremstyle{definition}
\newtheorem{definition}[theorem]{Definition}
\newtheorem{remark}[theorem]{Remark}

\newtheorem{example}[theorem]{Example}
\newtheorem{examples}[theorem]{Examples}
\newtheorem{exercise}[theorem]{Exercise}
\newtheorem{notation}[theorem]{Notation}
\def\b{\mathbb }

\def\phi{\varphi }

\def\scr{\mathscr}
\swapnumbers

\begin{document}
\title{Dunkl operators: Theory and applications}
\author{Margit R\"osler}

\address{Mathematisches Institut, Universit\"at G\"ottingen,
Bunsenstr.~3-5, D-37073 G\"ottingen, Germany}
\email{roesler@uni-math.gwdg.de}

\begin{abstract}

These lecture notes are intended as an introduction to the theory of
rational Dunkl operators and the associated special
functions, with an emphasis on positivity and asymptotics.
 We start with an outline of the general concepts: Dunkl operators,
the intertwining operator, the Dunkl kernel and the Dunkl transform.
We point out the connection with integrable particle systems
of Calogero-Moser-Sutherland type, and discuss some systems
of orthogonal polynomials associated with them. 
A major part is  devoted to positivity results for the intertwining
operator and the Dunkl kernel, the Dunkl-type heat semigroup, and 
related  probabilistic aspects.
The notes conclude with recent results on the asymptotics
of the Dunkl kernel.

\end{abstract}

\maketitle

\tableofcontents

%

\section{Introduction}
\label{S0}

While the theory of special functions in one variable has a long and rich history, 
the growing interest in  special functions of several variables is comparatively recent.
During the last years, there has in particular been a rapid 
development in the area of 
special functions with reflection symmetries and the harmonic analysis related with root systems. 
The motivation for this subject comes to 
some extent from the theory of  Riemannian symmetric spaces, 
whose spherical functions can be written as multi-variable special functions depending on certain
discrete sets of parameters.
A key tool in the study of special functions with reflection
symmetries are Dunkl operators. Generally speaking, these 
 are commuting 
differential-difference operators, associated to a finite 
reflection group on a Euclidean space. 
The first class of such operators, 
now often called ``rational'' Dunkl operators, 
were introduced 
by C.F. Dunkl in the late 80ies. In a series of papers ([D1-5]), he built up the
framework for a theory of special functions and integral transforms in several variables related
with reflection groups. Since then, 
various other classes of Dunkl operators have become important, in the first place the
trigonometric Dunkl operators of Heckman, Opdam and the Cherednik operators. These will not be discussed
in our notes; for an overview, we refer to \cite{He}.
An important  motivation 
 to study Dunkl operators originates  in 
their relevance for  the analysis of 
 quantum many body systems of  Calogero-Moser-Sutherland type. These describe algebraically integrable
systems  in one dimension and  have gained considerable  interest in 
mathematical physics, especially in conformal field theory. A good bibliography is contained in
\cite{DV}.

The aim of these lecture notes is an introduction to rational
Dunkl theory, with an emphasis on the author's results in this area. Rational Dunkl operators
bear a rich analytic structure which  is
 not only due to their commutativity, but
also to the existence of an intertwining operator between Dunkl operators  and
usual partial derivatives. We shall first give
an overview of the general concepts, including  an account on the relevance of Dunkl operators 
in the study of Calogero-Moser-Sutherland models. We
also  discuss some of the special functions related with them. 
A major topic will be positivity results; these concern  the  intertwining operator
as well as the kernel of the Dunkl transform, and lead to a variety of positive 
semigroups in the Dunkl setting with possible probabilistic interpretations.
 We make this  explicit at hand of the most 
important example: the Dunkl-type heat semigroup, which is generated by the analog of the 
Laplacian in the Dunkl setting. The last section presents recent results on
the asymptotics of the Dunkl kernel  and the short-time behavior 
of heat kernels associated with root systems.

\section{Dunkl operators and the Dunkl transform} 
\label{S1}

The aim of this section is to provide an introduction to
 the theory of rational Dunkl operators, which we shall call Dunkl operators for
short,  and to the Dunkl transform.   General references 
are [D1-5], \cite{DJO}, \cite{dJ1} and \cite{Op1}; for a background on reflection 
groups and root systems the reader is 
referred to \cite{Hu} and \cite{GB}.
We do not intend to give a complete survey, but rather focus on those aspects which will 
be important in the context of this lecture series.

\subsection{Root systems and reflection groups}

The basic ingredient in the theory of Dunkl operators are root systems and finite reflection
groups, acting on some Euclidean space  $(E, \langle\,.,.\,\rangle)$ of finite dimension $N$.
  We shall always assume that $E=\b R^N$ with the standard Euclidean
 scalar product $\,\langle x,y\rangle \,=\, \sum_{j=1}^N x_j y_j$.
For $\alpha\in \b R^N\setminus\{0\}$, we denote by 
$\sigma_\alpha $ the reflection in
the hyperplane $\langle\alpha\rangle^\perp$ orthogonal to $\alpha$, i.e.
\[\sigma_\alpha(x)\,=\, x -
2\,\frac{\langle\alpha,x\rangle}{|\alpha|^2}\,\alpha\,,\]
where $|x|:=\sqrt{\langle x,x\rangle}$. Each reflection $\sigma_\alpha$ is contained 
in the orthogonal group $O(N, \b R)$. 
We start with the basic definitions:

 \begin{definition}
Let $ R\subset \b R^N\setminus\{0\}$  be a finite set. Then $R$ is called a 
 \textit{root system}, if \parskip=-1pt
\begin{enumerate}
\item[\rm{(1)}] $R\cap \,\b R\alpha = \{\pm\alpha\}$ for all $\alpha\in R$;
\item[\rm{(2)}] $\sigma_\alpha(R) = R$ for all $\alpha\in R.$
\end{enumerate}
The subgroup $\, G= G(R)\subseteq O(N, \b R)$ which is generated by the 
reflections $\{\sigma_\alpha,\,\alpha\in R\}$ 
is called the \textit{reflection group}
(or \textit{Coxeter-group}) associated with $R$. The dimension of $span_{\b R}R$
is called the rank of $R$. 
\end{definition}

Property  (1) is called reducedness. It is often not required in Lie-theoretic contexts, where 
 instead the root systems under consideration are assumed to  be crystallographic. This means that 
\[ \frac{2\langle\alpha,\beta\rangle}{\langle\beta,\beta\rangle}\,\in \b Z \quad\text{for all }\, \alpha,\beta\in R.\]
If $R$ is crystallographic and has full rank, then 
$span_{\b Z} R$ forms a lattice in $\b R^N$ (called the \textit{root-lattice}) which is 
 stabilized by the action of the associated reflection group.

\begin{lemma}\begin{enumerate}\itemsep=-1pt
\item[\rm{(1)}] For any root system $R$ in $\b R^N,$ the reflection group $G=G(R)$ is finite.
\item[\rm{(2)}] The set of reflections contained in $G(R)$ is exactly 
$\{\sigma_\alpha\,,\,\alpha\in R\}.$
\end{enumerate}
\end{lemma}

\begin{proof} As $R$ is left invariant by $G$, we have a  natural homomorphism 
$\varphi: G \to S(R)$
of $G$ into the symmetric group of $R$, given by $\varphi(g)(\alpha):= g\alpha \in R$. 
This homomorphism is injective: indeed, each reflection $s_\alpha$, and therefore also
each element  $g\in G$ fixes pointwise the orthogonal complement of the subspace spanned by $R$.
If also $g(\alpha) = \alpha$ for all $\alpha\in R,$ then $g$ must be the identity. 
This implies assertion (1) because the order of $S(R)$ is finite. Property  (2) 
is more involved. An elegant proof can be found in Section 4.2 of \cite{DX}. 
\end{proof}

\begin{exercise} If $g\in O(N,\b R)$ and $\alpha\in \b R^N\setminus\{0\},$ then 
$g\,\sigma_\alpha g^{-1} = \sigma_{g\alpha}.$
\end{exercise}

Together with part (2) of the previous lemma, this shows that  there
is a bijective correspondence between  the conjugacy classes of reflections in $G$ and the
orbits in $R$ under the natural action of $G$. We shall need some more concepts: 
Each root system can  be written as a disjoint union $R= R_+ \cup (-R_+),$ 
where $R_+$ and $-R_+$ are separated by a hyperplane through the origin. Such a set $R_+$ is
called a \textit{positive subsystem}. Of course, its choice is not unique. 
The set of 
hyperplanes $\{ \langle\alpha\rangle^\perp, \,\alpha\in R\}$ divides $\b R^N$ into connected 
open components, called the \textit{Weyl chambers} of $R$. 
It can be shown that  the topological closure 
$\overline C$ of any chamber $C$ is a fundamental domain for $G$,
i.e. $\overline C$ is naturally homeomorphic with the 
 space $(\b R^N)^G$ of all $G$-orbits in $\b R^N$, endowed
with the quotient topology. $G$ permutes the reflecting hyperplanes as well as the chambers.

\begin{exercise} \textit{Dihedral groups.} 
In the Euclidean plane $\b R^2$, let $d\in O(2,\b R)$ denote the rotation
around $2\pi/n$ with $n\geq 3$  and $s$ the reflection at the $y$-axis. Show that 
the group $\mathcal D_n$ generated by $d$ and $s$ consists of all orthogonal transformations 
which preserve a regular $n$-sided polygon centered at the origin. (Hint: $dsd = s$.)
Show that $D_n$ is a finite reflection group and determine its root system. 
Can the crystallographic condition  always be satisfied?
\end{exercise}

\begin{examples} 
\textbf{(1) Type} $A_{N-1}$. Let $S_N$ denote the  symmetric group in $N$ elements. 
It acts faithfully 
 on $\b R^N$ by permuting the standard basis vectors $e_1,\ldots,e_N$. Each
transposition $(ij)$ acts as a reflection $\,\sigma_{ij}$ sending 
$ e_i - e_j$ to its negative. Since $S_N$ is generated by transpositions, it is a finite reflection group. 
A root system of $S_N$ is given by
\[ R = \{ \pm(e_i - e_j), \,1\leq i<j\leq N\}.\]
Its span is $(e_1 + \ldots e_N)^\perp$, and thus the rank is $N-1$.

\noindent
\textbf{(2) Type} $B_N$. Here $G$ is the reflection group in $\b R^N$ generated by the 
transpositions $\sigma_{ij}$ as above, as well as the sign changes  
$\sigma_i: e_i\mapsto -e_i\,,\> i=1,\ldots,N.$ 
The group of sign changes is isomorphic to $\b Z_2^N$, intersects
$S_N$ trivially and is normalized by $S_N$, so  $G$ is isomorphic with the semidirect product
$S_N\ltimes \b Z_2^N$. The corresponding root system has rank $N$; it is given by
\[ R\,=\{ \pm e_i,\, 1\leq i \leq N, \, \pm(e_i\pm e_j), 1\leq i < j\leq N\}.\]
\end{examples}

A root system $R$ is called \textit{irreducible}, if it cannot be written as the 
orthogonal disjoint union $R= R_1\cup R_2$ of two root systems $R_1\,,\, R_2$. Any root system can be uniquely written as an orthogonal disjoint union of irreducible root systems.
There exists a classification of all irreducible root systems in terms of 
Coxeter graphs. The crystallographic ones are made up by $4$ infinite families 
 $A_n,\, B_n$ (those  discussed above),   $C_n\,, D_n$, as well as  $5$ exceptional 
root systems.
For details, we refer to \cite{Hu}.  

\subsection{Dunkl operators}

Let $R$ be a fixed root system in $\b R^N$ and $G$ the associated reflection group. 
From now on we assume that  $R$ is normalized in the sense that 
$\langle \alpha,\alpha\rangle =2$ for all 
$\alpha\in R$; this simplifies formulas, but is no loss of generality
for our purposes. 
The Dunkl operators attached with $R$ can be considered as 
perturbations of the usual partial derivatives by reflection parts. 
These reflection parts are coupled by parameters, which  are given in terms of a 
multiplicity function:

\begin{definition} A function $k:R\to \b C$ on the root system $R$ is called a 
\textit{multiplicity function} on $R$, if it is invariant under the natural action of $G$ on $R$.
The $\b C$-vector space of multiplicity functions on $R$ is denoted by $K$.
\end{definition}

Notice that 
the  dimension of $K$ 
is equal to the number of $G$-orbits in $R$. We write $k\geq 0$ if
 $k(\alpha)\geq 0$ for all $\alpha\in R$.

\begin{definition} Let $k\in K.$ Then for $\xi\in \b R^N,$ 
the \textit{Dunkl operator} $T_\xi:= T_\xi(k)$ is defined (for $f\in C^1(\b R^N)$) by
\[ T_\xi f(x):= \partial_\xi f(x) + \sum_{\alpha\in R_+} 
   k(\alpha)\,\langle\alpha, \xi\rangle\, 
    \frac{f(x) - f(\sigma_\alpha x)}{\langle\alpha, x\rangle}.\]
Here $\partial_\xi$ denotes the directional  derivative corresponding
    to $\xi$, and $R_+$ is a fixed positive subsystem of $R$.
For the $i$-th standard basis vector $\xi=e_i\in \b R^N$ we use 
the abbreviation $T_i=T_{e_i}.$ 
\end{definition} 

The above definition does not depend on the special choice of $R_+$, 
thanks to the $G$-invariance of $k$. 
In case $k=0$, the $T_\xi$ reduce to the 
corresponding directional
derivatives. The operators $T_\xi$ were introduced and first studied for $k\geq 0$ by
C.F. Dunkl ([D1-5]). They enjoy regularity properties similar
to usual partial derivatives on various spaces of smooth functions on $\b R^N$. 
We shall use the following notations:

\begin{notation}
\begin{enumerate}\itemsep=1pt
\item[\rm{1.}] $\b Z_+ := \{0,1,2,\ldots\}$.
\item[\rm{2.}] $\Pi := \b C[\b R^N]$ is the $\b C$-algebra of polynomial
functions on $\b R^N$. It  has a natural grading 
\[ \Pi =
\bigoplus_{n\geq 0} \mathcal P_n\,,\]
 where $\mathcal P_n$ is the subspace of homogeneous polynomials of
(total) degree $n$. 
\item[\rm{3.}] $\scr S(\b R^N)$ denotes the Schwartz space of rapidly decreasing functions
on $\b R^N$, 
\[ \mathscr S(\b R^N):= \{ f\in C^\infty(\b R^N): \,
\|x^\beta\partial^\alpha\! f\|_{\infty, \b R^N}\, <\infty \quad\text{for all }\, \alpha, \,\beta\in \b Z_+^N\}.\]
It is a Fr\'echet space with the usual locally convex topology. 
\end{enumerate}
\end{notation}

\noindent
The Dunkl operators $T_\xi$ have the following regularity properties:

\begin{lemma}\label{L:regular} 
\begin{enumerate}
\item[\rm{(1)}] If $f \in C^m(\b R^N)$ with $m\geq 1$,  then $T_\xi f\in
  C^{m-1}(\b R^N)$.  
\item[\rm{(2)}] $T_\xi $ leaves $C_c^\infty(\b R^N)$ and $\scr S(\b R^N)$ invariant. 
\item[\rm{(3)}] $T_\xi $ is homogeneous of degree $-1$ on $\Pi$,
  that is, $T_\xi\,p\in \mathcal P_{n-1}$ for $p\in \mathcal
  P_n$. 
\end{enumerate}
\end{lemma}

\begin{proof} All statements follow from the representation 
\[ \frac{f(x)-f(\sigma_\alpha x)}{\langle\alpha,x\rangle}\,=\,
\int_0^1\partial_\alpha
f\big(x-t\langle\alpha,x\rangle\alpha\big)dt \quad\text{ for
  }\>  f\in C^1(\b R^N),
\> \alpha\in R\]
(recall our normalization $\langle\alpha,\alpha\rangle = 2$). 
(1) and (3) are immediate; the proof of 
(2) (for $\mathscr S(\b R^N)$) is also straightforward but more technical; it 
can be found in \cite{dJ1}.
\end{proof}


Due to the $G$-invariance of $k$, the Dunkl operators $T_\xi$ are $G$-equivariant:
In fact, consider the 
 natural action of  $O(N,\b R)$  on
functions $f:\b R^N\to \b C$, given by
\[ h\cdot f(x):= f(h^{-1}x), \quad h\in O(N, \b R).\]
Then an easy calculation shows: 

\begin{exercise}\label{E:equiv}
$\,\displaystyle g\circ T_\xi \circ g^{-1}\,=\, T_{g\xi}  \quad\text{for all }\,
g\in G.$
\end{exercise}

\noindent
Moreover, there holds a product rule:

\begin{exercise}
If $f,g\in C^1(\b R^N)$ and at least
one of them is $G$-invariant, then
\begin{equation}\label{(1.2)} 
 T_\xi(fg)\,=\, T_\xi(f)\cdot g \,+\, f\cdot T_\xi(g).
\end{equation}
\end{exercise}

The most striking property of the Dunkl operators, which is the foundation for rich analytic 
structures related with them, is the following

\begin{theorem}\label{T:commut} \enskip  For fixed $k,$ the associated
$T_\xi = T_\xi(k), \>\xi\in \b R^N\,$ commute. 
 \end{theorem}

This result was obtained in \cite{Du1} by a clever direct argumentation. 
An alternative proof, relying on Koszul complex ideas,  is given in \cite{DJO}.
As a consequence of Theorem \ref{T:commut} there exists an algebra homomorphism 
$\Phi_k: \Pi \,\to\,\text{End}_{\b C}(\Pi)$ which is defined by 
\[ \Phi_k: x_i\mapsto T_i, \,\, 1\mapsto id.\]
For $p\in \Pi$ we write
\[ p(T):= \Phi_k(p).\]

The classical case $k=0$ will be distinguished by the notation $\Phi_0(p) =: p(\partial)$. 
Of particular importance is the  $k$-\textit{Laplacian}, which is 
defined by 
\[ \Delta_k := p(T) \quad\text{with }\, p(x) = |x|^2.\]

\begin{theorem}\label{T:aplacedarst}
 \begin{equation}\label{(1.3)} \Delta_k \,=\, \Delta  + 2\sum_{\alpha\in R_+}
k(\alpha)\delta_\alpha \quad\,\text{ with }\quad \delta_\alpha f(x)\,=\,
\frac{\langle\nabla f(x),\alpha\rangle}{\langle\alpha ,x\rangle} - 
\frac{f(x)-f(\sigma_\alpha x)}{\langle\alpha , x\rangle^2};
\end{equation}
 here $\Delta$ and $\nabla$ denote the usual Laplacian and gradient 
respectively. 
\end{theorem}

\noindent
This representation is obtained by a direct calculation (recall again our convention 
$\langle\alpha,\alpha\rangle = 2$ for all $\alpha\in R$) by use of the following Lemma:

\begin{lemma}\cite{Du1}  For $\alpha\in R$, define 
\[\rho_\alpha f(x):= \frac{f(x) - f(\sigma_\alpha x)}{\langle\alpha,x\rangle} 
\quad (f\in C^1(\b R^N)).\]
Then
\[ \sum_{\alpha,\beta\in R_+} k(\alpha)k(\beta) \langle\alpha,\beta\rangle \rho_\alpha\rho_\beta
\,=\, 0.\]
\end{lemma}

\noindent
It is not difficult to check that 
\[ \Delta_k = \sum_{i=1}^N T_{\xi_i}^2
\]
for any orthonormal basis $\{\xi_1,\ldots,\xi_N\}$ of $\b R^N$, see 
\cite{Du1} for the proof.
Together with the 
 $G$-equivariance of the Dunkl operators, this immediately
implies that $\Delta_k$ is $G$-invariant, i.e.
\[g\circ\Delta_k = \Delta_k\circ g \quad (g\in
G).\]

\begin{examples} \label{E:dunklops}
 {\bf (1) The rank-one case.} \enskip In case $N=1$, the only 
choice of $R$ is $R=\{\pm \sqrt 2\}$, which is the root system of type $A_1$.
The corresponding  reflection group is
$G=\{id,\sigma\}$ acting on $\b R$ by $\sigma(x)=-x$. The
Dunkl operator
$T:= T_{1}$ associated with the multiplicity parameter 
$k\in \b C$ is given by
\[ Tf(x)\,=\, f^\prime(x) + k\,\frac{f(x)-f(-x)}{x}.\]
Its square $T^2$, when restricted to the even subspace 
$C^2(\b R)^e\, := \{f\in C^2(\b R): 
f(x)=f(-x)\}$ is given by a singular Sturm-Liouville operator:
\[ T^2\vert_{C^2(\b R)^e}f\,=\, f^{\prime\prime} + 
\frac{2k}{x}\cdot f^\prime \,.\]

{\bf (2) Dunkl operators  of type $A_{N-1}$.} Suppose $G=S_N$ with root system of type $A_{N-1}$.
(In contrast to the above example, $G$ now acts on $\b R^N$).
As all transpositions  are conjugate in $S_N$,  
 the vector space of multiplicity
functions  is one-dimensional. The Dunkl operators associated
with the multiplicity parameter $k\in \b C$ are given by
\[ T_i^S \,=\,\partial_i \,+\,k\cdot\sum_{j\not= i}
\frac{1-\sigma_{ij}}{x_i-x_j} \quad ( i=1,\ldots, N),\]
and the $k$-Laplacian is
\[ \Delta_k^S\,=\, \Delta +2k\sum_{1\leq i<j\leq N}
\frac{1}{x_i-x_j}\Big[(\partial_i-\partial_j) -
\frac{1-\sigma_{ij}}{x_i-x_j}\Big]\,.\]

{\bf (3)  Dunkl operators of type $B_N$.} Suppose $R$ is a root system of type $B_N$, 
corresponding to $G= S_N\ltimes \b Z_2^N$.  There are two conjugacy classes of
 reflections in $G$, leading to multiplicity functions of the form $k =
 (k_0, k_1)$ with $k_i\in \b C$. The associated Dunkl operators are given by
\[ T_i^B\,=\, \partial_i\, +\, k_1\frac{1-\sigma_i}{x_i}\,+\,
k_0\cdot\sum_{j\not=i} \Big[\frac{1-\sigma_{ij}}{x_i-x_j} + \frac{1-\tau_{ij}}{x_i+x_j}
\Big] \quad (i=1,\ldots, N),\]
where $\tau_{ij}:= \sigma_{ij}\sigma_i\sigma_j\,.$
\end{examples}

\subsection{A  formula of Macdonald and its analog in Dunkl theory}

In the classical theory of spherical harmonics (see for instance \cite{Hel}) the 
following bilinear pairing on $\Pi$, sometimes called Fischer product, plays an important role:
\[ [p,q]_0 := \bigl(p(\partial)q\bigr)(0), \quad p,q\in \Pi.\]
This pairing  is closely related to the scalar product 
in $L^2(\b R^N, e^{-|x|^2/2}dx)$; in fact, 
in his short note \cite{M2} Macdonald observed the following  identity:
\[ [p,q]_0\,=\, (2\pi)^{-N/2} \int_{\b R^N} e^{-\Delta/2} p(x)\, 
e^{-\Delta/2} q(x)\,e^{-|x|^2/2}\,dx.\]
Here $e^{-\Delta/2}$ is  well-defined as a linear operator
 on $\Pi$ by means of the  terminating series
\[ e^{-\Delta/2}p = \sum_{n=0}^\infty \frac{(-1)^n}{2^n n!}\Delta^n p.\]
Both the Fischer product as well as 
Macdonald's identity have a useful generalization in the Dunkl setting.
In the following, 
we shall always restrict to  the case 
$k\geq 0.$

\begin{definition} For $p, q\in \Pi$ define
\[ [p,q]_k\,:=\, \bigl(p(T)q\bigr)(0).\]
\end{definition}

\noindent
This bilinear form was introduced in \cite{Du2}.
We collect some of its basic properties:

\begin{lemma}\label{L:Skalprod}\parskip=-1pt
\begin{enumerate}
\item[\rm{(1)}] If $p\in \mathcal P_n$ and $q\in \mathcal P_m$ 
    with $n\not=m$, then $\,[p,q]_k =0.$
\item[\rm{(2)}] $\displaystyle [x_i\,p\,,q]_k\,=\, [p,T_i\,q]_k 
     \quad (p,q\in \Pi,\> i=1,\ldots, N).$
\item[\rm{(3)}] $\displaystyle [g\cdot p\,,g\cdot q]_k\,=\, [p,q]_k 
         \quad (p,q\in \Pi,\> g\in G).$
\end{enumerate}
\end{lemma}

\begin{proof} (1) follows from the homogeneity of the Dunkl operators, (2) is clear from the definition, and  (3) follows from Exercise \ref{E:equiv}. 
\end{proof}

Let $w_k$ denote the weight function on $\b R^N$ defined by 
\begin{equation}\label{(1.1b)}
 w_k(x)\,=\, \prod_{\alpha\in R_+} |\langle
\alpha,x\rangle|^{2k(\alpha)}.\end{equation}
It is $G$-invariant and homogeneous of degree $2\gamma$, with
the index
\begin{equation}\label{(1.1a)}
 \gamma:=\gamma(k):=\,\sum_{\alpha\in R_+} k(\alpha).
\end{equation}
Notice that by $G$-invariance of $k$, we have  $k(-\alpha) = 
k(\alpha)$ for all
 $\alpha\in R$. Hence  this definition does again not depend on the special
choice of  $R_+$. Further, we define the constant
\[ c_k:= \int_{\b R^N} e^{-|x|^2/2} w_k(x)dx,\]a
a so-called Macdonald-Mehta-Selberg integral. 
There exists a closed form for it which was conjectured and proved
by Macdonald \cite{M1} for the infinite series of root systems. 
An extension to arbitrary crystallographic reflection groups is due to Opdam \cite{Op1}, 
and there are computer-assisted proofs for some non-crystallographic root systems.
As far as we know, a general proof for arbitrary root systems has not yet been found. 

\noindent
We shall need the following anti-symmetry of the Dunkl operators:

\begin{proposition}\label{P:Antisym} \cite{Du3} \enskip Let $k\geq 0$. Then
  for every $f\in \scr S(\b R^N)$ and $g\in C_b^1(\b R^N)$,
\[ \int_{\b R^N} T_\xi f(x)g(x)w_k(x)dx\,=\, -\int_{\b R^N}
f(x)
T_\xi g(x)w_k(x)dx\,.\]
\end{proposition}

\begin{proof} A short calculation. In order to have the appearing
integrals well defined, one has to assume $k\geq 1$ first, and then extend the result
to general $k\geq 0$ by analytic continuation.
\end{proof}

\begin{proposition}\label{P:Macdonald}
For all $p,q\in \Pi$,
\begin{equation}\label{(1.10)}
 [p,q]_k\,=\, c_k^{-1} \int_{\b R^N} e^{-\Delta_k/2} p(x)\, 
e^{-\Delta_k/2} q(x)\,e^{-|x|^2/2}\,w_k(x)dx.
\end{equation}
\end{proposition}

\noindent
This result is due to Dunkl (\cite{Du2}). As the Dunkl Laplacian is homogeneous of degree $-2$, 
the operator $e^{-\Delta_k/2}$ is well-defined and bijective on $\Pi$, and it
 preserves the degree.
We give here 
 a direct proof which is partly taken from an unpublished part of M. de Jeu's
thesis (\cite{dJ3}, Chap. 3.3). It involves the following commutator
results in $\text{End}_{\b C}(\Pi)$, where as usual, $\,[A,B] = AB-BA\,$ for
$A, B\in \text{End}_{\b C}(\Pi)$.

\begin{lemma}\label{L:Commut} For $i=1,\ldots,N,$ \parskip=-1pt
\begin{enumerate}\itemsep=1pt
\item[\rm{(1)}] $\displaystyle
  \big[x_i\,,\,\Delta_k/2\big]\,=\, -T_i\,;$
\item[\rm{(2)}]  $\displaystyle \big[x_i\,,\,e^{-\Delta_k/2}\big]\,=\, T_i\,
  e^{-\Delta_k/2}$.
\end{enumerate}
\end{lemma}

\begin{proof} (1) follows by direct calculation, c.f. \cite{Du1}. Induction then yields that
\[ \big[x_i\,,(\Delta_k/2)^n\big]\,=\,
-n\,T_i\bigl(\Delta_k/2\bigr)^{n-1} \quad \text{for }n\geq 1,\]
and this implies (2).
\end{proof}

\begin{proof}[Proof of Proposition \ref{P:Macdonald}] 
Let $i\in \{1,\ldots, N\},$ and denote the right side of \eqref{(1.10)}
 by $(p,q)_k$. 
Then by the anti-symmetry of $T_i$ in
$L^2(\b R^N, w_k),$ the product rule for $T_i$ and the above Lemma,
\begin{align}
(p,T_i\,q)_k\,=&\,\, c_k^{-1}\int_{\b R^N}
e^{-\Delta_k/2} p \cdot
\bigl(T_i\,e^{-\Delta_k/2}
q\bigr)\,e^{-|x|^2/2}\,w_k\,dx \notag\\
=&\,\, -c_k^{-1}\int_{\b R^N} T_i\bigl(e^{-|x|^2/2}
e^{-\Delta_k/2}p\bigr)\cdot \bigl(e^{-\Delta_k/2}q\bigr)\, w_k\,dx\notag\\
=&\,\, c_k^{-1} \int_{\b R^N} e^{-\Delta_k/2}(x_i\,p)\cdot
\bigl(e^{-\Delta_k/2}q\bigr)\, e^{-|x|^2/2}\,w_k\,dx\,=\,
(x_i\,p,q)_k.\notag
\end{align}
But  the form $\,[\,.,.\,]_k\,$ has the same property by Lemma
\ref{L:Skalprod}(2). It is now easily checked that the assertion is true if $p$ or $q$ is constant,
and then, by induction on $\max(\text{deg}\,p,\text{deg}\,q)$, for all homogeneous $p, q$. 
This suffices by the linearity of both forms. \end{proof}

\begin{corollary} \label{C:Skalarprod} Let again  $k\geq 0$. Then 
the pairing  $\,[\,.\,,\,.\,]_k\,$ on $\Pi$ is symmetric and non-degenerate, i.e. 
$\, [p,q]_k = 0$ for all $q\in \Pi$ implies that $p=0$. 
\end{corollary}

\begin{exercise} Check the details in the proofs of Prop. \ref{P:Macdonald} and Cor. 
\ref{C:Skalarprod}.
\end{exercise}

\subsection{Dunkl's intertwining operator}

It was first shown in \cite{Du2} that for non-negative multiplicity functions, 
the associated commutative algebra of Dunkl operators is intertwined with the algebra of 
usual partial differential operators by a unique linear and homogeneous isomorphism on 
polynomials. A  thorough analysis in \cite{DJO} subsequently revealed that for general $k,$ 
such  an intertwining operator exists if and only if the common kernel of the $T_\xi$, 
considered as linear operators on $\Pi$, contains no "singular" polynomials
 besides the
constants. More precisely, the following characterization holds:

\begin{theorem}\label{T:Inter}\cite{DJO} \enskip  
 Let $\,K^{reg}:= \big\{k\in K: \bigcap_{\xi\in \b R^N}
  \text{Ker}\,T_\xi(k)\, = \b C\cdot 1\big\}\,.$  Then the following statements are equivalent:
\begin{enumerate}
\item[\rm{(1)}] $k\in K^{reg};$
\item[\rm{(2)}] There exists a unique linear isomorphism
  ("\textit{intertwining operator}") $V_k$ of
  $\Pi$ such that 
  \[ V_k(\mathcal
  P_n) = \mathcal P_n\,,\>\> V_k\,\vert_{ \mathcal P_0}\, =\, id
  \quad\text{ and }\>\> 
  T_\xi V_k\,=\, V_k\,\partial_\xi \quad\text{ for all }\,\xi\in \b
  R^N.\]
\end{enumerate}
\end{theorem}

The proof of this result is by induction on the degree of homogeneity and
requires only linear algebra.

\noindent
The intertwining operator $V_k$ commutes with the $G$-action:

\begin{exercise}
$\, \displaystyle g^{-1}\circ V_k\circ g\,=\, V_k \quad (g\in G).$

\noindent
\textit{Hint:} Use the $G$-equivariance of the $T_{\xi}$ and the defining properties of $V_k$.
\end{exercise}

\begin{proposition}
$\,\,\displaystyle \{k\in K: k\geq 0\}\subseteq K^{reg}.$
\end{proposition}

\begin{proof}
Suppose that $p\in \oplus_{n\geq 1}\mathcal P_n$ 
satisfies $T_\xi(k)p = 0$ for all $\xi\in \b R^N.$ Then
$[q,p]_k = 0$ for all $q\in \oplus_{n\geq 1}\mathcal P_n$, and hence also
$[q,p]_k = 0$ for all $q\in \Pi$. Thus $p=0$, by the
 non-degeneracy of $[\,.\,,\,.\,]_k$. 
\end{proof}

\noindent
The complete singular parameter set $ \,K\setminus K^{reg}$ is explicitly determined
 in \cite{DJO}. It is an open subset of $K$ which
 is invariant under complex conjugation, and contains 
$\{k\in K: \text{Re }k\geq 0\}.$
Later in these lectures,  we will in fact restrict our attention 
to non-negative multiplicity functions. These are of particular interest concerning 
our subsequent positivity results, which could not be expected for non-positive
multiplicities.
Though the intertwining operator plays an important role in Dunkl's theory, an 
explicit ``closed''  form for it is known so far only in some special cases. 
 Among these are
\begin{enumerate} 
\item[\rm{1.}] \emph{The rank-one case.}  Here 
\[K^{reg}\,=\, \b C\setminus\big\{ -1/2-n,\> n\in \b Z_+\big\}.\]
  The associated intertwining
  operator is given explicitly by
\[ V_k\bigl(x^{2n}\bigr)\,=\,
\frac{\bigl(\frac{1}{2}\bigr)_n}{\bigl(k+\frac{1}{2}\bigr)_n}\,
x^{2n}\,;\quad V_k\bigl(x^{2n+1}\bigr)\,=\,
\frac{\bigl(\frac{1}{2}\bigr)_{n+1}}{\bigl(k+\frac{1}{2}\bigr)_{n+1}}\,
x^{2n+1}\,,\]
where $\,(a)_n = \,\Gamma(a+n)/\Gamma(a)\,$ is the Pochhammer-symbol.
For $\text{Re}\,k >0$, this amounts to  the following integral representation 
(see \cite{Du2}, Th. 5.1):
\begin{equation}\label{(1.4)}
 V_k\, p(x)\,=\, \frac{\Gamma(k+1/2)}{\Gamma(1/2)\,\Gamma(k)}
  \int_{-1}^1 p(xt)\,(1-t)^{k-1}(1+t)^k\,dt.
\end{equation} 
\item[\rm{2.}] \emph{The case $G=S_3$.}
 This was studied in \cite{Du4}. Here 
\[K^{reg}= 
\b C\setminus\{-1/2-n, -1/3-n, -2/3-n, \, n\in
  \b Z_+\}.\]
\end{enumerate} 

In order to bring  $V_k$ into action in a further development of the theory, it is important
to extend it to larger function spaces. For this we shall always assume that  $k\geq 0$.
In a first step, $V_k$ is extended to a bounded linear operator on suitably normed algebras of homogeneous series 
on a ball. This concept goes back to \cite{Du2}.

\begin{definition}\label{D:ALg} For $r>0,$ let $B_r:=\{x\in \b R^N: |x|\leq r\}$ denote the closed
ball of
radius $r,$ and let $A_r$  be the closure of $\Pi$ with respect to the norm 
\[ \|p\|_{A_r}:= \sum_{n=0}^\infty \|p_n\|_{\infty, B_r} \, \quad\text{for }\, 
   p=\sum_{n=0}^\infty p_n\,, \,p_n\in \mathcal P_n.\]
\end{definition}

\noindent Clearly $A_r$ is a commutative Banach-$*$-algebra under the pointwise operations and with
complex conjugation as involution. 
Each $f\in A_r$ has a unique representation 
$f = \sum_{n=0}^\infty f_n\,$ with $f_n\in \mathcal P_n$, and is continuous on the ball $B_r$
and real-analytic in its interior. 
The topology of $A_r$ is stronger than the topology induced by the uniform norm on $B_r$. 
Notice also that
$A_r$ is not 
closed with respect to
$\|.\|_{\infty, B_r}$ and that $A_r\subseteq A_s$ with $\|\,.\,\|_{A_r} \geq \|\,.\,\|_{A_s}$
for $s\leq r$.

\begin{theorem}\label{T:Homog}
 $\, \|V_k p\|_{\infty, B_r} \leq \|p\|_{\infty, B_r}$ for each $p\in \mathcal P_n$. 
\end{theorem}

\noindent
The proof of this result is given in \cite{Du2} and can also be found in
 \cite{DX}.
 It
 uses the \textit{van der Corput-Schaake inequality}
which states that for each \textit{real-valued} $p\in \mathcal P_n$,
\[ \text{sup}\,\{ |\langle \nabla p(x),y\rangle|: \,x,y\in B_1\}\,
        \leq \,n \|p\|_{\infty, B_1}.\]
Notice that here the converse inequality is trivially satisfied, because 
$\,\langle \nabla p(x), x\rangle = np(x)$ for $p\in \mathcal P_n$. 
The following is now immediate:

\begin{corollary} $\|V_k f\|_{A_r} \,\leq \|f\|_{A_r}$ for every $f\in \Pi$, and 
$V_k$ extends uniquely to a bounded linear operator on $A_r$
 via 
\[ V_k f:= \sum_{n=0}^\infty V_k f_n \quad\text{for }\, f=\sum_{n=0}^\infty f_n\,.\]
\end{corollary}

Formula \eqref{(1.4)} shows in particular that in the rank-one case with $k>0,$
the operator $V_k$ is positivity-preserving on polynomials. It was conjectured by 
Dunkl in \cite{Du2} that for arbitrary reflection groups and non-negative multiplicity functions, 
the linear functional $f\mapsto V_kf(x)$ on $A_r$ should be positive. 
We shall see in Section 4.1. that this is in fact true.
As a consequence, we shall obtain the existence of 
a positive integral representation generalizing \eqref{(1.4)}, which in turn 
allows to extend $V_k$  to larger function spaces. This positivity result
 also has important consequences for the 
structure of the 
Dunkl kernel, which generalizes the usual exponential function in the Dunkl 
setting. We shall introduce it in the
following section.

\begin{exercise}\label{E:Banach}
The \textit{symmetric spectrum} $\Delta_S(A)$ of a (unital) 
commutative Banach-$*$-algebra $A$
is defined as the set  of all non-zero algebra 
homomorphisms $\varphi:A\to \b C$ satisfying the $*$-condition $\, \varphi(a^*) = 
\overline{\varphi(a)}$ for all $a\in A.$ It is a compact Hausdorff space with the 
weak-*-topology (sometimes called the Gelfand topology).

\smallskip
\noindent
Prove that the symmetric spectrum of the algebra $A_r$ is given by
\[ \Delta_S(A_r) = \{\varphi_x: \,x\in B_r\},\]
where $\varphi_x$ is the evaluation homomorphism $\varphi_x(f):= f(x)$. 
Show also that the mapping $x\mapsto \varphi_x$ is a homeomorphism from $B_r$ onto $\Delta_S(A_r)$. 
\end{exercise}

\subsection{The Dunkl kernel}

Throughout this section we assume that $k\geq 0$. Moreover, we denote by
$\langle\,.\,,\,.\,\rangle$ not only the Euclidean scalar product on $\b R^N,$ but also its
 \textit{bilinear} extension to  $\b C^N\times\b C^N.$ 

\noindent
For fixed $y\in \b C^N,$ the exponential function 
$x\mapsto e^{\langle x,y\rangle}$ belongs to each of the algebras $A_r\,,\, r>0.$
This justifies the following 

\begin{definition} \cite{Du2} For $y\in \b C^N$, define 
\[  E_k(x,y):= V_k\bigl(e^{\langle\,.\,,y\rangle})(x), \quad x\in \b R^N.\]
\end{definition}

\noindent
The function $E_k$ is called the \textit{Dunkl-kernel}, or 
$k$- exponential kernel, associated with $G$ and $k$. 
It can alternatively be characterized as the solution of a joint
eigenvalue problem for the associated Dunkl operators.

\begin{proposition}\label{P:Uniq} Let $k\geq 0$ and $y\in\b C^N.$ 
Then $f= E_k(\,.\,,y)$ is the unique solution of the
system 
\begin{equation}\label{(1.106)}
  T_\xi\,f\,=\, \langle \xi,y\rangle f \quad \text{for all }\xi\in \b R^N
\end{equation}
which is real-analytic on $\b R^N$ and satisfies $\,f(0)=1.$  
\end{proposition}

\begin{proof} $E_k(\,.\,y)$ is real-analytic on $\b R^N$ by our construction.
Define 
\[ E_k^{(n)}(x,y):= \frac{1}{n!}V_k\langle\,.\,,y\rangle^n(x), \quad x\in \b R^N, \,n=0,1,2,\ldots.\]
Then $\, E_k(x,y) = \sum_{n=0}^\infty   E_k^{(n)}(x,y),\,$ 
and the series converges uniformly and absolutely with respect to $x$.
The homogeneity of $V_k$ immediately implies $E_k(0,y)=1$. Further, by the intertwining property,
\begin{equation}\label{(1.103)}
 T_\xi E_k^{(n)}(\,.\,,y) \,=\, \frac{1}{n!} V_k\,\partial_\xi \langle\,.\,,y\rangle^n\,=\,
\langle \xi,y\rangle E_k^{(n-1)}(\,.\,,y)
\end{equation}
for all $n\geq 1$. This shows that $E_k(\,.\,,y)$ solves  \eqref{(1.106)}.
To prove uniqueness,  suppose that $f$ is 
a real-analytic solution of \eqref{(1.106)} with $f(0)=1$.
Then $T_\xi$ can be applied termwise to the homogeneous expansion $\, f= \sum_{n=0}^\infty f_n\,,
f_n\in \mathcal P_n$, and comparison of homogeneous parts shows that
\[ f_0 = 1, \quad T_\xi f_n = \langle\xi,y\rangle f_{n-1} \quad \text{for }\, n\geq 1.\]
As $\{ k\in K: k\geq 0\}\subseteq K^{reg}$, it follows by induction that all $f_n$ are 
uniquely determined. 

\end{proof}

While this construction has been carried out only for $k\geq 0$, there is a 
more general result by Opdam which assures the existence of a general exponential kernel
with properties according to the above lemma for arbitrary regular 
 multiplicity parameters.  The following is a weakened version of \cite{Op1}, Prop. 6.7;
it in particular implies that $E_k$ has a holomorphic extension to $\b C^N\times \b C^N$:

\begin{theorem}\label{T:Opdam}
  For each $k\in K^{reg}$ and $y\in \b C^N$, the system
\[ T_\xi\,f\,=\, \langle \xi,y\rangle f \quad (\xi\in \b R^N)\]
has a unique solution $\, x\mapsto E_k(x,y)\,$ 
which is real-analytic on $\b R^N$ and satisfies
$\,f(0)=1.$  
Moreover, the mapping $(x,k, y)\mapsto E_k(x,y)$ 
 extends to a meromorphic function on 
$\b C^N\times K\times\b C^N$ with pole set $\b C^N\times (K\setminus K^{reg})\times\b C^N$
\end{theorem}

\noindent
We collect some further  properties of the Dunkl kernel $E_k$. 

\begin{proposition}\label{P:Kernequiv} Let $k\geq 0,\, x,y\in \b C^N, \,\lambda\in \b C$ and $g\in G$. 
\begin{enumerate}
\item[\rm{(1)}] $E_k(x,y) = E_k(y,x)$ 
\item[\rm{(2)}] $E_k(\lambda x,y) = E_k(x,\lambda y)\, $ and $\,E_k(gx,gy) = E_k(x,y).$
\item[\rm{(3)}] $\overline{E_k(x,y)} = E_k(\overline x,\overline y).$ 
\end{enumerate}
\end{proposition}

\begin{proof} (1) This is shown in \cite{Du2}. (2) is easily obtained from the definition
of $E_k$ together with the homogeneity and equivariance properties of $V_k$. 
For (3), notice that $f:= \overline{E_k(\,.\,,y)},$ which is again real-analytic on $\b R^N,$ 
satisfies $\,T_\xi f = 
\langle\xi,\overline y\rangle\,f, \, f(0) =1$. By the uniqueness part of the above Proposition,
 $\overline{E_k(x,y)} = E_k(x,\overline y)$ for all real $x$. Now both 
$x\mapsto \overline{E_k(\overline x,y)}$ and $x\mapsto 
E_k(x,\overline y)$ are holomorphic on $\b C^N$ and agree on $\b R^N$. Hence they coincide.
\end{proof}

\noindent
Just as with the intertwining operator, the kernel $E_k$ is explicitly
known  for some particular cases only.  An important example is again the rank-one 
situation:

\begin{example}\label{E:1kern}  In the \textit{rank-one case}  
with $\text{Re}\,k >0,$  
the integral representation
  \eqref{(1.4)} for $V_k$ implies that for all $x,y\in \b C$, 
\[ E_k(x,y)\,=\, \frac{\Gamma(k+1/2)}{\Gamma(1/2)\,\Gamma(k)}\int_{-1}^1
e^{txy}(1-t)^{k-1}(1+t)^k\,dt\, =\, e^{xy}\cdot
\phantom{}_1F_1(k,2k+1,-2xy).\]
This can also be written as
\[ E_k(x,y)\,=\, j_{k-1/2}(ixy)\,+\,\frac{xy}{2k+1}\,
j_{k+1/2}(ixy), \]
where for $\alpha \geq -1/2$, $j_\alpha$ is the normalized spherical
Bessel function
\begin{equation}\label{(1.8)}
 j_\alpha(z)\,=\,
2^\alpha\Gamma(\alpha+1)\cdot\frac{J_\alpha(z)}{z^\alpha}\,=\,
\Gamma(\alpha+1)\cdot\sum_{n=0}^\infty
\frac{(-1)^n(z/2)^{2n}}{n!\,\Gamma(n+\alpha+1)}\,.\end{equation}
\end{example}

\noindent
This  motivates the following

\begin{definition} \cite{Op1}
The $k$-\textit{Bessel function} (or generalized Bessel function) is defined by
\begin{equation}\label{(5.0)}
 J_k(x,y) := \frac{1}{|G|} \sum_{g\in G} E_k(gx,y) \quad (x,y\in \b C^N).
\end{equation} 
\end{definition}

Thanks to Prop. \ref{P:Kernequiv} $J_k$ is $G$-invariant in both arguments and
therefore naturally considered on Weyl chambers of $G$ 
(or their complexifications). 
In the rank-one case, we  have
\[ J_k(x,y) = j_{k-1/2}(ixy).\]

\noindent
It is a well-known fact from classical analysis that 
for fixed $y\in \b C$, the function $f(x) = j_{k-1/2}(ixy)$ is the unique
analytic solution of the differential equation 
\[ f^{\prime\prime} + \frac{2k}{x} f^\prime \,=\, y^2 f \]
which is even and normalized by $f(0)=1$. In order to see how this can be generalized 
to the multivariable case, 
consider the algebra of $G$-invariant polynomials on $\b R^N$,
\[ \Pi^G = \{ p\in \Pi: g\cdot p = p \quad\text{for all }\, g\in G\}.\]
If $p\in \Pi^G$, then it follows from the equivariance of the 
Dunkl operators (Exercise \ref{E:equiv}) 
that $p(T)$ commutes with the $G$-action;  a 
detailed argument for this
is given in \cite{He2}. Thus $p(T)$ leaves $\Pi^G$ invariant, and we obtain in particular that
 for fixed $y\in \b C^N$, the $k$-Bessel function 
$J_k(\,.\,,y)$ is a solution to the following
Bessel-system: 
\[ p(T)f\,=\, p(y)f \quad \text{for all }\, p\in \Pi^G, \, \, f(0)=1.\]
According to \cite{Op1}, it is in fact the only  $G$-invariant and analytic solution.
We mention that there exists a  group theoretic context in which, for a certain parameters $k$, 
generalized Bessel functions 
occur in a natural way: namely  as the spherical functions of a Euclidean type symmetric space, associated with a so-called Cartan motion group. We refer
to  \cite{Op1} for this
connection and to \cite{Hel} for the necessary background in semisimple Lie theory.

The Dunkl kernel is of particular interest as it gives rise to an associated
integral transform on $\b R^N$ which generalizes the Euclidean 
Fourier transform in a natural way. This transform will be discussed in the
following section. 
Its  definition and essential properties  rely on
suitable growth estimates for  $E_k$. In our case $k\geq 0$, 
the best ones to be expected are available:

\begin{proposition}\label{P:Roughestim}\cite{R3}
 For all $x\in \b R^N, \, y\in \b C^N$ and all multi-indices $\alpha\in \b Z_+^N$, 
\[\vert \partial_y^\alpha E_k(x,y)\vert\,\leq\, |x|^{|\alpha|}\max_{g\in G} 
     e^{\text{Re}\langle gx,y\rangle}. \]
In particular,
$\,\,\,\vert E_k(-ix,y)\vert\,\leq\, 1 \>\>$ for all $x,y\in \b R^N$. 
\end{proposition}

\noindent
This result will be obtained later from a positive integral representation of
Bochner-type for $E_k$, c.f. Cor. \ref{C:Bochner}.
 M. de Jeu had slightly weaker bounds in \cite{dJ1},
 differing by
an additional factor $\sqrt{|G|}.$

\noindent
We conclude this section with two important reproducing properties for the Dunkl kernel.
Notice that the above estimate on $E_k$ assures the
convergence of the involved integrals.

\begin{proposition}\label{P:Reprod} \enskip Let $k\geq
  0.$  Then
\begin{enumerate}\itemsep=1pt
\item[\rm{(1)}] 
$\,\displaystyle 
\int_{\b R^N} e^{-\Delta_k/2}p\,(x)\, E_k(x,y)\,
 e^{-|x|^2/2}w_k(x)dx\,=\,c_k\, e^{\langle y,y\rangle/2} p(y) 
     \quad (p\in \Pi,\> y\in \b C^N).$
\item[\rm{(2)}] $\displaystyle \, 
\int_{\b R^N} E_k(x,y)\,E_k(x,z)\,e^{-|x|^2/2}w_k(x)dx\,=\, c_k\,
e^{(\langle y,y\rangle + \langle z,z\rangle)/2} E_k(y,z) \quad (y,z\in
\b C^N).$
\end{enumerate}
\end{proposition}

\begin{proof} (c.f. \cite{Du3}.) We shall use the Macdonald-type formula 
\eqref{(1.10)} for the pairing $[\,.\,,\,.\,]_k$. First, we prove that 
\begin{equation}\label{(1.104)}
 [E_k^{(n)}(x,\,.\,), \,.\,]_k \,=\, p(x) 
\quad\text{for all }\, p\in \mathcal P_n, \, x\in \b R^N.
\end{equation}
In fact, if $p\in \mathcal P_n,$ then 
\[p(x) = (\langle x, \partial_y\rangle^n/n!)p(y) \,\quad\text{and }\,\, 
V_k^x p(x) = E_k^{(n)}(x, \partial_y)p(y).\] 
Here the uppercase index in $V_k^x$ denotes the relevant variable.
Application of $V_k^y$ to both sides 
 yields $\,V_k^x p(x) = \, E_k^{(n)}(x,T^y) V_k^y p(y).$
As $V_k$ is bijective on $\mathcal P_n$, this implies \eqref{(1.104)}.
For fixed $y,$ let $L_n(x):= \sum_{j=0}^n E_k^{(j)}(x,y).$ If $n$ is larger than the degree
of $p,$ it follows from \eqref{(1.104)} that
$\,[L_n, p]_k \,=\, p(y).$ Thus in view of the Macdonald formula, 
\[c_k^{-1}\int_{\b R^N} e^{-\Delta_k/2} L_n(x) e^{-\Delta_k/2}p(x)
 e^{-|x|^2/2}w_k(x)dx\,=\,p(y).\]  
On the other hand, it is easily checked that
\[ \lim_{n\to\infty} e^{-\Delta_k/2} L_n(x) \,=\, e^{-\langle y,y\rangle/2}E_k(x,y).\] 
This gives (1). 
Identity (2) then follows from (1), again by homogeneous expansion of $E_k$. 
\end{proof}

\subsection{The Dunkl transform}

The Dunkl transform was introduced in 
\cite{Du3} for non-negative multiplicity functions 
and further studied in \cite{dJ1} in the more general case
$\text{Re}\,k\geq 0$. In these notes, we again restrict ourselves to $k\geq 0$.

\begin{definition} The Dunkl transform associated
with $G$ and $k\geq 0$ is given by
\[ \widehat .^{\,k}\,:\, L^1(\b R^N, w_k) \,\to\, C_b(\b R^N); \quad 
\widehat f^{\,k}(\xi):=  c_k^{-1}\int_{\b R^N} f(x) E_k(-i\xi,x) 
\,w_k(x)dx \>\>(\xi\in \b R^N).\]
The inverse transform is defined by $\,f^{\vee k}(\xi)=\widehat f^{\,k}(-\xi)$.
\end{definition}

Notice that  $\widehat f^{\, k}\in C_b(\b R^N)$ results from our
 bounds on  $E_k$.
The Dunkl transform shares many properties with the classical Fourier transform.
Here are the most basic ones:

\begin{lemma} Let $f\in \mathscr S(\b R^N).$ Then for $j=1,\ldots, N,$ 
\begin{enumerate}\itemsep=1pt
\item[\rm{(1)}] $\displaystyle \widehat f^{\,k} \in C^\infty(\b R^N)$ and
$\displaystyle  T_j(\widehat f^{\,k})\,=\, -(ix_j f)^{\wedge k}.$
\item[\rm{(2)}]
 $\displaystyle (T_j f)^{\wedge k}(\xi) = i\xi_j \widehat f^{\,k}(\xi).$
\item[\rm{(3)}] The Dunkl transform leaves $\scr S(\b R^N)$ invariant. 
\end{enumerate}
\end{lemma}

\begin{proof} (1) is obvious from \eqref{(1.106)}, and
(2) follows from the anti-symmetry relation  (Prop. \ref{P:Antisym}) for the Dunkl operators.
For (3), notice that it suffices to prove that 
$\partial^\alpha_\xi (\xi^\beta \widehat f^{\,k}(\xi))$ is bounded for arbitrary multi-indices 
$\alpha,\,\beta.$ By the previous Lemma, we have $\xi^\beta \widehat f^{\,k}(\xi)= \widehat g^{\,k}(\xi)$ for some $g\in \scr S(\b R^N).$ Using the  growth 
bounds of Proposition \ref{P:Roughestim}  yields the assertion.

\end{proof}

\begin{exercise} \begin{enumerate}
\item[\rm{(1)}] $C_c^\infty(\b R^N)$ and  $\scr S(\b R^N)$ are dense in 
$L^p(\b R^N, w_k), \, \, p=1,2.$
\item[\rm{(2)}] Conclude the 
\textit{Lemma of Riemann-Lebesgue} for the Dunkl transform: 
\[ f\in L^1(\b R^N, w_k) \,\Longrightarrow\, \widehat f^{\,k} \in C_0(\b R^N).\]
Here $C_0(\b R^N)$ denotes the space of continuous functions on $\b R^N$ which vanish at
infinity. 
\end{enumerate}
\end{exercise}

\noindent
The following are the main results for the Dunkl transform; we omit the proofs but refer the
reader to   \cite{Du3} and  \cite{dJ1}:

\begin{theorem}\label{P:Dunkltrafo}
\begin{enumerate}
\item[\rm{(1)}]
The Dunkl transform $\,f\mapsto \widehat f^{\,k}$ is a
  homeomorphism of 
$\scr S(\b R^N)$ with period $4$.
\item[\rm{(2)}] (Plancherel theorem) The Dunkl transform has a unique
  extension to an isometric 
isomorphism of $L^2(\b R^N, w_k)$. We denote this isomorphism again by  
$\,f\mapsto \widehat f^{\,k}$. 
\item[\rm{(3)}]  ($L^1$-inversion) For all $\,f\in L^1(\b R^N, w_k)$ with 
   $\,\widehat f^{\,k}\in L^1(\b R^N, w_k),$
\[ f\,=\, (\widehat f^{\,k}\,)^{\vee k} \quad \text{a.e.}\]
\end{enumerate}
\end{theorem}

\section{CMS models and generalized Hermite polynomials}

\subsection{Quantum Calogero-Moser-Sutherland models}

Quantum Calogero-Moser-Sutherland (CMS) models describe quantum mechanical
 systems of $N$ identical particles on a circle or line which interact pairwise
 through long range potentials of inverse square type. They are exactly solvable
 and have gained considerable interest in theoretical physics during 
the last years. Among the broad literature in this area, we refer to \cite{DV}, 
  \cite{LV}, \cite{K},
\cite{BHKV}, \cite{BF1}-\cite{BF3}, \cite{Pa}, \cite{Pe}, \cite{UW}, \cite{Du7}.
CMS models  have in particular attracted some attention  
in conformal field theory, and they are being used to 
test the ideas of fractional statistics (\cite{Ha}, \cite{Hal}).
While explicit spectral
resolutions of  such models were already obtained by Calogero and
 Sutherland (\cite{Ca}, \cite{Su}), a new aspect in the understanding of
 their algebraic structure and quantum integrability 
was much later initiated by
 \cite{Po} and \cite{He2}. The Hamiltonian under consideration is
 hereby modified by 
 certain exchange operators, which allow to write  it in a
 decoupled form. These exchange modifications can be expressed in terms
 of  Dunkl operators of type $A_{N-1}$. The  Hamiltonian of the
\textit{linear CMS model with harmonic confinement} in $L^2(\b
R^N)$ is given by
\begin{equation}\label{(3.70)}
 \mathcal H_C\,=\, -\Delta \,+\,g\sum_{1\leq i<j\leq N} 
\frac{1}{(x_i-x_j)^2}\, +\omega^2|x|^2; \end{equation}
here $\omega>0$ is a frequency parameter  and  $g\geq -1/2$ is
a coupling constant. In case $\omega =0$, \eqref{(3.70)} describes the
free Calogero model. On the other hand, if $g=0$, then $\mathcal H_C$ coincides
with the Hamiltonian  of the $N$-dimensional isotropic harmonic oscillator, 
\[ \mathcal H_0 = -\Delta + \omega^2 |x|^2.\]
The spectral decomposition of this operator in $L^2(\b R^N)$ is well-known: 
The spectrum is discrete, $\sigma(\mathcal H_0) =\{(2n+N)\omega,\, n\in \b Z_+\}$, 
and the  classical multivariable Hermite functions (tensor products 
of one-variable Hermite functions, c.f. Examples \ref{E:Ex1}), form a complete set of eigenfunctions. 
The study of the Hamiltonian $\mathcal H_C$  was initiated by  Calogero  (\cite{Ca}); he
 computed its spectrum and determined  the structure of the 
bosonic eigenfunctions and scattering states in the confined and free case, respectively.
Perelomov \cite{Pe} observed that
\eqref{(3.70)} is completely
quantum integrable, i.e. there exist $N$ commuting, algebraically independent
symmetric linear operators in $L^2(\b R^N)$ including $\mathcal H_C$.
We mention that the complete 
integrability of  the
classical Hamiltonian systems associated with \eqref{(3.70)} 
 goes back to Moser \cite{Mo}.
There exist generalizations of the classical Calogero-Moser-Sutherland
models in the context of
abstract root systems, see for instance \cite{OP1},  \cite{OP2}.  In particular, if
$R$ is an arbitrary root system on $\b R^N$  and $k$ is a nonnegative multiplicity
function on it, then the  corresponding  abstract  Calogero
Hamiltonian  with harmonic confinement is  given by 
\[ \widetilde{\mathcal H}_k\,=\, -\widetilde {\mathcal F}_k + \omega^2|x|^2\]
with the formal expression
\[ \widetilde{\mathcal F}_k\,=\, \Delta - \,2\sum_{\alpha\in R_+}
k(\alpha)(k(\alpha)-1) \frac{1}{\langle \alpha,x\rangle^2}.\] 
If $R$ is of type $A_{N-1}$, then $\widetilde{\mathcal H}_k$  just coincides with
$\mathcal H_C$. For both the classical and the
quantum case, partial
results on the integrability of this model are due to 
Olshanetsky and Perelomov \cite{OP1}, \cite{OP2}.
A  new aspect in the understanding of the  algebraic structure
 and the  quantum
integrability of CMS systems was  initiated by 
Polychronakos \cite{Po} and  Heckman \cite{He2}. The underlying  idea 
is to construct quantum integrals
for  CMS models from differential-reflection  operators. 
Polychronakos introduced 
them in terms of an  ``exchange-operator formalism''
for \eqref{(3.70)}.  He thus  obtained a complete
 set of commuting observables for
\eqref{(3.70)} in an elegant way.
In \cite{He2} it was observed in
 general that
 the complete algebra of quantum integrals 
for free, abstract Calogero models is intimately connected with 
the corresponding algebra
 of Dunkl operators. 
Let us briefly describe this connection: 
Consider the following modification of 
$\widetilde{\mathcal F}_k$, involving reflection terms:
\begin{equation}\label{(3.71)}  \mathcal F_k\,=\, \Delta -2\sum_{\alpha\in R_+}
 \frac{k(\alpha)}{\langle\alpha,x\rangle^2}\,(k(\alpha)
 -\sigma_\alpha)\,.\end{equation}
In order to avoid singularities
 in the reflecting hyperplanes, it is suitable to carry out a gauge
 transform by $w_k^{1/2}.$ A short calculation, using again
results from \cite{Du1}, gives
\[ w_k^{-1/2} \mathcal F_k w_k^{1/2} \,=\, \Delta_k\,,\]
c.f. \cite{R4}.
Here $\Delta_k$ is the Dunkl Laplacian associated with $G$ and $k$. 
Now consider the algebra of $\Pi^G$ of $G$-invariant polynomials on $\b R^N$.  By a 
classical theorem of Chevalley (see e.g. \cite{Hu}), it is generated by $N$ 
homogeneous, algebraically independent elements. 
For $p\in \Pi^G$ we denote by $\text{Res}\,(p(T))$ the restriction of the Dunkl operator
$p(T)$ to $\Pi^G$ (Recall that $p(T)$ leaves  $\Pi^G$ invariant!). 
Then  
\[ \mathcal A:= \bigl\{ \text{Res}\,p(T): \,p\in \Pi^G\}\]
is a commutative algebra of differential operators on $\Pi^G$ containing 
the operator
\[ \text{Res}\,(\Delta_k)\,=\, w_k^{-1/2}\widetilde{\mathcal F}_k\, w_k^{1/2},\]
and $\mathcal A$ has $N$ algebraically independent generators, called quantum integrals
for the free Hamiltonian  $\widetilde{\mathcal F}_k$.
\subsection{Spectral analysis of abstract CMS Hamiltonians}

This section is devoted to a spectral analysis of abstract linear CMS operators with harmonic confinement. We follow the expositions in \cite{R2}, \cite{R5}.
To simplify formulas, we fix $\omega=1/2$;  
corresponding results for general $\omega$ can always be obtained by rescaling. 
We again work with the gauge-transformed version with reflection terms, 
\[\mathcal H_k:= \,w_k^{-1/2} (-\mathcal F_k+ \frac{1}{4}|x|^2)\,w_k^{1/2}\,=\, 
   -\Delta_k + \frac{1}{4}|x|^2.\]
Due to the anti-symmetry of the first order Dunkl operators (Prop. \ref{P:Antisym}), 
   this operator is symmetric and densely defined 
in $L^2(\b R^N, w_k)$ with domain $\mathcal D(\mathcal H_k):= \scr S(\b
R^N)$. Notice that in case $k=0$, $\mathcal H_k$  is just the 
 Hamiltonian of the $N$-dimensional isotropic harmonic oscillator. We further 
consider the Hilbert space $L^2(\b R^N, m_k)$, where $m_k$ is 
 the probability measure 
\begin{equation}\label{(3.90)}
 dm_k\,:= \,c_k^{-1} e^{-|x|^2/2}w_k(x)dx
\end{equation}
and the operator
\[ \mathcal J_k := \, -\Delta_k + \sum_{i=1}^N x_i\partial_i \]
in $L^2(\b R^N, m_k)$, with domain 
$\mathcal D(\mathcal J_k):= \Pi$. It can be shown by standard methods that $\Pi$ is dense
in $L^2(\b R^N, m_k)$. We do not carry this out; a proof can be found in \cite{R3} or in
\cite{dJ4}, where a comprehensive treatment of density questions in several variables is given.

The next theorem contains a  complete description of the spectral properties 
of $\mathcal H_k$ and $\mathcal J_k$ and 
generalizes the already mentioned well-known 
facts for the classical harmonic oscillator Hamilonian.
For the proof, we shall employ the    
$sl(2)$-commutation relations of the operators
\[ E:= \frac{1}{2}|x|^2,\>\> F:= -\frac{1}{2}\Delta_k \quad\text{and}\>\> 
H:= \sum_{i=1}^N x_i\partial_i +\left(\gamma +N/2\right)\]
on $\Pi$ (with the index $\gamma=\gamma(k)$ as defined in
\eqref{(1.1a)}) which can be found in    \cite{He2}. They are 
\begin{equation}\label{(3.5)}
\big[H,E\big] = 2E,\>\> \big[H,F\big] = -2F,\>\>\big[E,F\big] = H.
\end{equation}
Notice that the first two relations are immediate consequences of the fact that the 
Euler operator  
\begin{equation}\label{(3.5a)}\rho:=\sum_{i=1}^N x_i\partial_i\,
\end{equation}
 satisfies  
$\,\rho(p)=np\,$ for each homogeneous $p\in \mathcal P_n$. We start with the following

\begin{lemma} On $\,\mathcal D(\mathcal J_k) = \Pi,$  
\[ \mathcal J_k\,=\, e^{|x|^2/4}(\mathcal H_k - (\gamma+N/2))
\,e^{-|x|^2/4}.\] 
In particular, $\mathcal J_k$ is symmetric in 
$L^2(\b R^N, m_k)$.
\end{lemma}

\begin{proof} From \eqref{(3.5)} it is easily verified by induction that 
\[ \big[\Delta_k, E^n\big] = 2nE^{n-1}H + 2n(n-1)
E^{n-1}\quad\text{for all }
\>\>n\in \b N,\]
and therefore $ \,\big[\Delta_k, e^{-E/2}\,
\big] = -e^{-E/2}H + \frac{1}{2} E e^{-E/2}.\,$
Thus on $\Pi$, 
\[ \mathcal H_k\, e^{-E/2} \, = \,-\Delta_k e^{-E/2}
\, + \frac{1}{2} Ee^{-E/2}  = 
\,-e^{-E/2}\Delta_k\,  +   e^{-E/2} H \,=\, e^{-E/2} 
(\mathcal J_k + \gamma+N/2).\]
\end{proof}

\begin{theorem} \label{T:Spectral}
 The spaces $L^2(\b R^N, m_k)$ and $L^2(\b R^N, w_k)$
admit  orthogonal Hilbert space decompositions into eigenspaces of the operators $\mathcal J_k$ and $\mathcal H_k$ respectively. More precisely, define
\[V_n\,:= \{e^{-\Delta_k/2}\, p\,: p\in \mathcal P_n\,\}
\subset \Pi, \quad
 W_n\,:= \{ e^{-|x|^2/4}\,q(x), \>\> q\in V_n\}\,\subset \scr
 S(\b R^N).\]
Then $V_n$ is the eigenspace of $\mathcal J_k$ corresponding to the 
eigenvalue $n$,  $W_n$ is the eigenspace of $\mathcal H_k$ 
 corresponding to the eigenvalue $\,n + \gamma+N/2, $ and 
\[ L^2(\b R^N, m_k)\,=\, \bigoplus_{n\in \b Z_+} V_n\,, \quad 
   L^2(\b R^N, w_k)\,=\, \bigoplus_{n\in \b Z_+} W_n\,.  \]
\end{theorem}

\begin{remark} 
A densely defined linear operator $(A, \mathcal D(A))$ in a Hilbert space $H$ 
is called essentially self-adjoint, if it satisfies
\begin{enumerate}
\item[\rm{(i)}] $A$ is symmetric, i.e. $\langle Ax,y\rangle = \langle x,Ay\rangle$
for all $x\in \mathcal D(A);$
\item[\rm{(ii)}] The closure $\overline A$ of $A$ is selfadjoint.
\end{enumerate}
In fact, every symmetric operator $A$ in $H$ has a unique closure $\overline A$ (because 
$A\subseteq A^*$, and the adjoint $A^*$ is closed).
If $H$ has a countable orthonormal basis $\{v_n,\,n\in \b Z_+\} \subset \mathcal D(A)$ consisting
of eigenvectors of $A$ corresponding to eigenvalues $\lambda_n\in \b R$, 
 then it is straightforward that $A$ is essentially
self-adjoint, and that the spectrum of the self-adjoint operator 
$\overline A$ is given by $\sigma(\overline A)= \{\lambda_n\,,n\in \b Z_+\}$. (See for instance
 Lemma
1.2.2 of \cite{Da}).

In our situation, the operator $\mathcal H_k$ is densely defined and symmetric in 
$L^2(\b R^N, w_k)$ (the first order Dunkl operators being anti-symmetric), and the
same holds for $\mathcal J_k$ in $L^2(\b R^N, m_k)$. The above theorem implies that 
$\mathcal H_k$ and $\mathcal J_k$ are essentially self-adjoint and that
\[ \sigma(\overline{\mathcal H_k}) = \{ n+\gamma+N/2, n\in \b Z_+\}, \quad
\sigma(\overline{\mathcal J_k}) = \b Z_+\,. \]
\end{remark}

\begin{proof}[Proof of Theorem \ref{T:Spectral}]  Equation
\eqref{(3.5)} and induction yield the commuting relations 
$\, \big[\rho,\Delta_k^n\big] = -2n\Delta_k^n\, $ for all $
n\in \b Z_+$, and hence   
\[\big[\rho, e^{-\Delta_k/2}\big]\, = \,\Delta_k\, e^{-\Delta_k/2}.\] 
If  $q\in \Pi$ is arbitrary and $\,p:= e^{\Delta_k/2}q\,,$ it  follows 
that 
\[\rho(q) \,=\, (\rho e^{-\Delta_k/2})(p)\, =\, e^{-\Delta_k/2}\rho(p) + 
\Delta_k e^{-\Delta_k/2\,}p\,=\, 
e^{-\Delta_k/2}\rho(p) + \Delta_k q.\]
Hence for $a\in \b C$ there are equivalent:
\[ (-\Delta_k + \rho)(q) = a q \,\Longleftrightarrow \, \rho(p) = ap\,
\Longleftrightarrow\, a=n\in \b Z_+ \>\>\text{and}\>\> p\in \mathcal P_n.\]
Thus each function from $V_n$ is an eigenfunction of $\mathcal J_k$ corresponding to 
the eigenvalue $n,$ and $\, V_n \perp V_m$ for 
$n\not=m$ by the symmetry of $\mathcal J_k$. This proves the statements for $\mathcal J_k$
because $\,\Pi  \,=\,\bigoplus V_n\,$ is dense 
in $L^2(\b R^N, m_k)$. 
 The statements for $\mathcal H_k$ are 
then  immediate by  the previous Lemma. 
\end{proof}

\subsection{Generalized Hermite polynomials}

The eigenvalues of the CMS Hamiltonians $\mathcal H_k$ and $\mathcal J_k$ are
highly degenerate if $N>1$. In this section, we construct natural orthogonal bases for
them. They are made up by generalizations of the
 classical $N$-variable Hermite polynomials and Hermite functions
to the 
Dunkl setting. We follow \cite{R2}, but change our normalization by a factor $2$.

\noindent
The starting point for our construction 
is the Macdonald-type  identity: if $p, q\in \Pi$, then
\begin{equation}\label{(3.110)}
 [p,q]_k\,=\, \int_{\b R^N} e^{-\Delta_k/2}p(x) e^{-\Delta_k/2}q(x) dm_k(x),
\end{equation}
with the probability measure $m_k$ defined according to \eqref{(3.90)}. 
Notice that $[\,.\,,\,.\,]_k$ is a scalar product on the $\b R$- vector space $\Pi_{\b R}$ of 
 polynomials with real coefficients.
Let $\{\phi_\nu\,, \nu\in \b Z_+^N\}$ be an orthonormal basis of $\Pi_{\b R}$ 
with respect to the scalar product $[.\,,.]_k$ such that $\phi_\nu\in 
\mathcal P_{|\nu|}$. As homogeneous polynomials of different degrees are orthogonal,
 the $\phi_\nu$ with fixed $|\nu|=n$ can for example be 
constructed by Gram-Schmidt orthogonalization within $\mathcal P_n\cap \Pi_{\b R}$  from an 
arbitrary ordered real-coefficient basis. If $k=0$, the
 canonical choice of the basis $\{\phi_\nu\}$ is just $\,\phi_\nu(x):= 
(\nu!)^{-1/2} x^\nu.$

\begin{definition}
The generalized Hermite polynomials $\{H_\nu\,, \>\nu\in \b Z_+^N\}$  
associated with the basis $\{\phi_\nu\}$ on $\b R^N$ are given by 
\begin{equation}
H_\nu(x):= e^{-\Delta_k/2}\phi_\nu(x).
\end{equation}
Moreover, we define the generalized Hermite functions on $\b R^N$ by
\begin{equation}
h_\nu(x):= e^{-|x|^2/4}H_\nu(x), \quad \nu\in \b Z_+^N.
\end{equation}
\end{definition}

\medskip
$H_\nu$ is a polynomial of degree $|\nu|$ satisfying $H_\nu(-x) = 
(-1)^{|\nu|}H_\nu(x)$ for all $x\in \b R^N$.  By virtue  of 
\eqref{(3.110)}, the $H_\nu, \>\nu\in \b Z_+^N$ form an 
orthonormal basis of $L^2(\b R^N, m_k).$

\begin{examples}\label{E:Ex1} \begin{enumerate}\itemsep=1pt
\item[\rm{(1)}]  \textit{Classical multivariable Hermite polynomials.}\enskip Let 
  $k=0$, and choose the standard orthonormal system $\phi_\nu(x) = (\nu!)^{-1/2}x^\nu$,
with respect to $[\,.\,,.\,]_0$. The associated Hermite
  polynomials are given by
\begin{equation}\label{(3.91)}
 H_\nu(x)\,=\, \frac{1}{\sqrt{\nu!}}\prod_{i=1}^N
  e^{-\partial_i^2/2}(x_i^{\nu_i})\,=\,
  \frac{2^{-|\nu|/2}}{\sqrt{\nu!}}\,\prod_{i=1}^N \widehat
  H_{\nu_i}(x_i/\sqrt 2),\end{equation}
where the $\widehat H_n, \> n\in \b Z_+$ are the classical Hermite
  polynomials on $\b R$ defined by 
\[ \widehat H_n(x)\,=\, (-1)^n\,e^{x^2}\,\frac{d^n}{dx^n}e^{-x^2}\,.\]
\item[\rm{(2)}] \textit{The one-dimensional case.} Up to sign changes, there exists only one  
    orthonormal basis
   with respect to $[\,.\,,\,.\,]_k$. The
  associated Hermite polynomials are given, up to multiplicative constants, by  
  the generalized Hermite polynomials $H_n^k(x/\sqrt{2})$  on $\b R$. These polynomials 
can be found  in \cite{Chi} and were further studied in \cite{Ros} in connection with a 
Bose-like oscillator calculus. The $H_n^k$
  are orthogonal with respect to $|x|^{2k} e^{-|x|^2}$ and can be
  written as
\[ \begin{cases}
     H_{2n}^k(x) = (-1)^n 2^{2n} n! \,L_n^{k-1/2}(x^2),\\
     H_{2n+1}^k(x) = (-1)^n 2^{2n+1} n! \,xL_n^{k+1/2}(x^2);
   \end{cases}\]
here  the $L_n^\alpha$ are the usual Laguerre polynomials of index 
$\alpha\geq -1/2$,  given by
\[ L_n^\alpha(x) = \frac{1}{n!} x^{-\alpha} e^x\,\frac{d^n}{dx^n}
\Bigl( x^{n+\alpha} e^{-x}\Bigr).\]
\item[\rm{(3)}] \textit{The $A_{N-1}$-case.} 
There exists a natural orthogonal system $\{\varphi_\nu\}$, made up by the so-called
\textit{non-symmetric Jack polynomials.} 
For a  multiplicity parameter $k> 0$, the associated non-symmetric 
Jack polynomials $E_\nu\,,\> \nu\in \b Z_+^N$, as introduced in \cite{Op} 
(see also \cite{KS}),
 are uniquely defined
by the following conditions:\parskip=-1pt
\begin{enumerate}
\item[\rm{(i)}] $\displaystyle E_\nu(x)\,=\, x^\nu \,+\, \sum_{\mu <_P\,\nu} 
c_{\nu,\,\mu} x^\mu\,$ with $c_{\nu,\mu}\in \b R$;
\item[\rm{(ii)}] For all $\mu <_P\,\nu$, $\displaystyle 
\bigl(E_\nu(x), x^\mu\bigr)_k\,= 0\,$
\end{enumerate}
\smallskip\noindent
Here $\, <_P\, $ is a dominance order defined within multi-indices 
of equal total length (see \cite{Op}), and 
the inner product $(.,.)_k$ on $\Pi\cap \Pi_{\b R}$ is given by
\[ (f,g)_k\,:=\, \int_{\b T^N} f(z) g(\overline z)\, \prod_{i<j} |z_i-z_j|^{2k} dz\]
with $\b T=\{z\in \b C: |z|=1\}$ and $dz$ being the Haar measure on $\b T^N$. 
If $f$ and $g$ have different total degrees, then $(f,g)_k=0$. The set $\{E_\nu\,,\> |\nu|=n\}$ forms a vector space basis of 
$\mathcal P_n \cap \Pi_{\b R}$. It can be shown (by use of $A_{N-1}$-type Cherednik 
operators) that the Jack polynomials $E_\nu$  are also 
orthogonal with respect to the Dunkl pairing 
$[\,.\,,\,.\,]_k$; for details see \cite{R2}. The corresponding generalized Hermite polynomials and their 
symmetric counterparts  
have been 
studied in \cite{L1}, \cite{L2} and in \cite{BF1} - \cite{BF3}.
\end{enumerate}
\end{examples}

As an immediate consequence of Theorem \ref{T:Spectral} we obtain analogues of the classical second order differential equations for generalized Hermite polynomials and Hermite functions:

\begin{corollary} \begin{enumerate}
\item[\rm{(i)}] $\displaystyle \bigl(-\Delta_k + \sum_{i=1}^N x_i\partial_i\bigr) H_\nu \,=\, |\nu| H_\nu\,.$
\item[\rm{(ii)}] $\displaystyle \bigl(-\Delta_k +\frac{1}{4}|x|^2\bigr) h_\nu\,=\, (|\nu|+\gamma+N/2) h_\nu\,.$
\end{enumerate}
\end{corollary}

Various further useful properties of the classical Hermite polynomials and Hermite functions have
extensions to our general setting. We conclude this section with a list of them. The proofs 
can be found  in  \cite{R2}. For further results on generalized Hermite polynomials, one
can also see for instance \cite{vD}.

\begin{theorem} Let $\{H_\nu\}$  be the Hermite polynomials and Hermite 
functions associated with the basis $\{\phi_\nu\}$ on $\b R^N$ and let  $x,y\in \b R^N$. Then 
\begin{enumerate}\itemsep=1pt
\item[\rm{(1)}]  $\displaystyle H_\nu(x) \,=\, (-1)^{|\nu|} e^{|x|^2/2}
\varphi_{\nu}(T) e^{-|x|^2/2}\,$  (Rodrigues-Formula)
\item[\rm{(2)}] $\displaystyle e^{-|y|^2/2} E_k(x,y)\,=\, \sum_{\nu\in \b Z_+^N} H_\nu(x)\varphi_\nu(y)\, $ (Generating relation)
\item[\rm{(3)}] (Mehler formula) For all $0<r<1,$ 
\[\sum_{\nu\in \b Z_+^N} H_\nu(x) H_\nu(y) r^{|\nu|}\,=\, \frac{1}{(1-r^2)^{\gamma +N/2}}\,
{\rm exp}\left\{-\frac{r^2(|x|^2 + |y|^2)}{2(1-r^2)}\right\} 
E_k\!\Bigl(\frac{rx}{1-r^2},\,y\Bigr).\]
\end{enumerate}
The sums on the left are absolutely convergent in both cases. 
\end{theorem}
The Dunkl kernel $E_k$ in (2) and (3) replaces the usual exponential function.
It comes in via the following relation with the (arbitrary!) basis $\{\varphi_\nu\}$:
\[ E_k(x,y) \,=\, \sum_{\nu\in \b Z_+^N} \varphi_\nu(x)\varphi_\nu(y) \quad (x,y\in \b R^N).\]

\begin{proposition}
The generalized Hermite functions $\{h_\nu\,,\, \nu\in \b Z_+^N\}$
are a basis of eigenfunctions of the Dunkl transform on $L^2(\b R^N, w_k)$ with
\[ h_\nu^{\wedge k} \,=\, (-i)^{|\nu|} h_\nu.\]
\end{proposition}

\section{Positivity results}

\subsection{Positivity of Dunkl's intertwining operator}

In this section it is always assumed that $k\geq 0.$ The reference is \cite{R3}.

We shall say that
a linear operator $A$ on $\Pi$ is  \textit{positive}, if $A$ leaves the positive cone 
\[\Pi_+ :=\{p\in \Pi: p(x)\geq 0 \quad \text{for all }\,x\in \b R^N\}\] 
invariant. 
The following theorem is the central result of this section:

\begin{theorem}\label{T:Pos1} $V_k$ is positive on $\Pi$. 
\end{theorem}

Once this is known, 
more detailed information about $V_k$ can be obtained by its extension  to the
 algebras $A_r$, which were introduced in Definition \ref{D:ALg}. This leads to

\begin{theorem} \label{T:Main2}
For each $x\in \b R^N$ there exists a unique probability measure
$\mu_x^k$ on
the Borel-$\sigma$-algebra of $\b R^N$ such that 
\begin{equation}\label{(4.11)}
 V_kf(x)\,=\, \int_{\b R^N} f(\xi)\,d\mu_x^k(\xi)\quad \text{ for all}\>\, 
f\in A_{|x|}.
\end{equation} 
The representing measures $\mu_x^k$ are compactly
supported with 
$\,{\rm supp}\,\mu_x^k \subseteq \, {\rm co}\,\{gx,\>g\in G\},$ the convex
hull of the orbit of $x$ under $G$. 
Moreover, they satisfy
\begin{equation}\label{(4.12a)}
\mu_{rx}^k(B)\,=\,\mu_x^k(r^{-1}B),\quad \mu_{gx}^k(B)\,=\,\mu_x^k(g^{-1}(B))
\end{equation} 
for each $r>0, \> g\in G$  and each Borel set $B\subseteq \b
R^N$.
\end{theorem}

The proof of Theorem \ref{T:Pos1} affords several steps, 
the crucial one being a reduction from the $N$-dimensional to a one-dimensional problem.
We shall give  an outline, but beforehand we turn to the
proof of Theorem \ref{T:Main2}. 

\begin{proof}[{\bf Proof of Theorem \ref{T:Main2}}]   
Fix  $x\in \b R^N$ and put $r=|x|.$ Then the mapping 
\[ \Phi_x : f\mapsto V_kf(x)\]
is a bounded linear functional on $A_r$, and  Theorem \ref{T:Pos1} implies that
it is positive on the dense subalgebra $\Pi$ of $A_r$, i.e. $\Phi_x(|p|^2) \geq 0.$ 
Consequently, 
$\Phi_x$ is a  positive functional on the full Banach-$*$-algebra
$A_r$. There exists a representation theorem of Bochner for 
 positive functionals on commutative 
Banach-$*$-algebras (see for instance Theorem 21.2 of \cite{FD}). It implies in our case
that  there exists a
unique  measure $\nu_x\in M_b^+(\Delta_S(A_r))$ such that 
\[
\Phi_x(f)\,=\, \int_{\Delta_S(A_r)} \widehat
f(\phi)\,d\nu_x(\phi)\quad\>\text{ for all }\> f\in A_r,
\]
with $\widehat f$  the Gelfand transform of $f$. Keeping Exercise \ref{E:Banach} in mind, one obtains 
 representing measures $\mu_x^k$ supported in the ball $B_{r}$; the sharper statement on the support is obtained by results of \cite{dJ1}. 
 The remaining statements are easy. 
\end{proof}

The key for the proof of Theorem \ref{T:Pos1} is a characterization of 
positive semigroups on polynomials which are generated by degree-lowering operators.
 We call a linear operator $A$ on $\Pi$ 
\textit{degree-lowering}, if $\,deg(Ap) < deg(p)$ for all $\Pi$. Again,  the exponential
$e^A\in \text{End}(\Pi^N)$ is defined by a terminating power-series, and it can be considered
as a linear operator on each of the finite dimensional spaces $\{p\in \Pi: \, deg(p)\leq m\}$.
Important examples of degree-lowering operators are linear operators which 
are homogeneous of some degree $-n <0$, such as Dunkl operators. 
The following key result characterizes positive semigroups generated 
by degree-lowering operators; it is an adaption of a well-known Hille-Yosida 
type 
characterization theorem for so called Feller-Markov semigroups which will be
discussed
a little later in our course, see Theorem \ref{T:Feller}.

\begin{theorem}\label{T:Posmin}
Let $A$ be a degree-lowering linear operator on $\Pi$. Then the following statements are 
 equivalent: \parskip=-2pt
\begin{enumerate}
\item[\rm{(1)}] $\, e^{tA}$ is positive on $\Pi$ for all $t\geq 0$.
\item[\rm{(2)}] $A$ satisfies the ``positive minimum principle''
\parskip=-1pt
\begin{center}
 {\rm (M)} \hspace{2pt} For every $p\in \Pi_+$ and $x_0\in \b R^N,$  
\hspace{2pt} 
  $p(x_0)=0\>$ implies $\, Ap\,(x_0)\geq 0.$ 
\end{center}
\end{enumerate} 
\end{theorem}

\begin{exercise}\label{E:lowering} 
\begin{enumerate}
\item[\rm{(1)}] 
Prove  implication \,(1) $\Rightarrow$ (2) of this theorem. 
\item[\rm{(2)}] Verify that the (usual) Laplacian $\Delta$ satisfies the positive minimum principle (M).
Can you extend this result to the Dunkl Laplacian $\Delta_k$? (C.f. Exercise \ref{E:mini}!)
\end{enumerate}
\end{exercise}

Let us now outline the proof of Theorem \ref{T:Pos1}. 
We consider the generalized Laplacian $\Delta_k$ associated
with $G$ and $k$, which is  homogeneous of degree $-2$ on $\Pi. $
With the notation 
introduced in \eqref{(1.3)}, it can be written as 
\begin{equation}\label{(4.5)}
 \Delta_k\,=\,\Delta\,+\, L_k\quad \text{ with }\> L_k\,=\,2\sum_{\alpha\in R_+} 
k(\alpha)\delta_\alpha\,.
\end{equation} 
Here 
 $\delta_\alpha$ acts in direction $\alpha$ only.

\begin{theorem}\label{H1} The operator $\, e^{-\Delta/2}e^{\Delta_k/2}$ is positive on $\Pi$.
\end{theorem}
 \begin{proof}
We shall deduce this statement from a positivity result for a suitable semigroup. For this, 
we  employ  Trotter's product formula, which works for degree-lowering operators just as
on finite-dimensional vector spaces:
 If $A, B$ are degree-lowering linear operators on $\Pi$, then 
\[e^{A+B} p\,(x)\,=\, \lim_{n\to\infty}
 \bigl(e^{A/n} e^{B/n}\bigr)^n p\,(x).\]
Thus, we can write 
\begin{align}
 e^{-\Delta/2} e^{\Delta_k/2}\, p\,(x)\,=\,&\, e^{-\Delta/2} e^{\Delta/2 + 
   L_k/2}\, p\,(x)\,= \,
  \lim_{n\to\infty}\,e^{-\Delta/2}\Bigl( e^{\Delta/2n}\, 
   e^{L_k/2n}\Bigr)^n p\,(x)\, \notag\\
 =\,& \lim_{n\to\infty} \prod_{j=1}^n \Bigl(e^{-(1-j/n)\cdot\Delta/2}\,
    e^{L_k/2n}\,e^{(1-j/n)\cdot\Delta/2}\Bigr) p\,(x).\notag 
 \end{align}
It therefore suffices to verify that the operators 
\[ e^{-s\Delta} e^{tL_k} e^{s\Delta}\,\quad (s,t\geq 0)\]
are positive on $\Pi.$ Consider $s$ fixed, then
\[ e^{-s\Delta} e^{tL_k} e^{s\Delta}\,=\, e^{tA} \quad\text{with }\, 
A = e^{-s\Delta}L_k e^{s\Delta}\,.\]
It is easily checked that $A$ is degree-lowering. Hence, in view of Theorem \ref{T:Posmin},
it remains to show that $A$ satisfies the positive minimum principle $(M)$. 
We may write
\[ A = \, e^{-s\Delta}L_k e^{s\Delta}\,=\, 2\sum_{\alpha\in R_+} k(\alpha) e^{-s\partial_\alpha^2}
\delta_\alpha\, e^{s\partial_\alpha^2};\]
here it was used that  $\delta_\alpha$ acts in direction $\alpha$ only. 
It can now be checked by direct computation that the one-dimensional operators 
$ e^{-s\partial_\alpha^2}\delta_\alpha e^{s\partial_\alpha^2}$ satisfy $(M),$ and as the $k(\alpha)$
are non-negative, this must be true for $A$ as well.
\end{proof}

\begin{proof}[{\bf Proof of Theorem \ref{T:Pos1}}]   
Notice first that  \begin{equation}\label{(3.31)} 
[V_k\, p,q]_k\,=\, [p,q]_0 \quad\text{ for all }\,p,q\in
\Pi.
\end{equation}
In fact, for $p,q\in \mathcal P_n$ with  $n\in \b Z_+$, one obtains
\[
[V_k\, p,q]_k\,=\, [q,V_k\, p]_k\,=\,q(T)(V_k\, p)\,=\, 
   V_k(q(\partial)p)\,=\, q(\partial)(p)\,=\,
[p,q]_0\,; \]
here the characterizing properties of $V_k$ and the fact that
   $q(\partial)(p)$ is a constant have been used. 
For general $p,q\in \Pi$, \eqref{(3.31)} then follows from the
orthogonality of the spaces $\mathcal P_n$ with
respect to both pairings.

Combining the Macdonald-type identity  \eqref{(1.10)} with part
\eqref{(3.31)}, we obtain for all $p,q\in
\Pi$ the identity
\[
c_k^{-1}\int_{\b R^N} e^{-\Delta_k/2}(V_k p)e^{-\Delta_k/2} q\,
e^{-|x|^2/2}\,w_k(x) dx\,=\, c_0^{-1}\int_{\b R^N} e^{-\Delta/2}
p \,e^{-\Delta/2} q \,e^{-|x|^2/2} dx.
\]
As $\displaystyle\, e^{-\Delta_k/2}(V_k\, p)\,=\,
V_k\bigl(e^{-\Delta/2}p\bigr)\,$, and as we may also replace $p$ by
$e^{\Delta/2}p\,$ and $q$ by $e^{\Delta_k/2}q\,$ in the above identity, it
follows that for all $p,q\in \Pi$
\begin{equation}\label{(2.20)} c_k^{-1}\int_{\b R^N} V_k p\,q\,e^{-|x|^2/2}\,w_k(x) dx\,=\, c_0^{-1}
\int_{\b R^N} p\, e^{-\Delta/2} e^{\Delta_k/2} q\,
e^{-|x|^2/2} dx.
\end{equation}
Due to Theorem \ref{H1}, the right  side of \eqref{(2.20)} 
is non-negative
for all $\,p,q\in \Pi_+$.  From this, the assertion can be deduced
 by standard density arguments ($\Pi$ is dense in $L^2(\b R^N, e^{-|x|^2/4}w_k(x)dx)$).
\end{proof}

\begin{corollary}\label{C:Bochner}
For each $y\in \b C^N$, the function $\,x\mapsto E_k(x,y)\,$ has the 
Bochner-type representation
\begin{equation}\label{(5.1)}
 E_k(x,y)\,=\,\int_{\b R^N} e^{\langle\xi,\,y\rangle} d\mu_x^k(\xi),
\end{equation}
where the $\mu_x^k$ are  the representing measures from Theorem \ref{T:Main2}. In
particular, $E_k$ satisfies the estimates stated in  Prop. \ref{P:Roughestim}, and 
\[E_k(x,y)>0 \quad\text{for all } x,y\in \b R^N.\]
Analogous statements hold for the $k$-Bessel function $J_k$. 
\end{corollary}

In those cases where 
 the generalized Bessel functions $J_k(.,y)$ allow an interpretation as the spherical functions 
of a Cartan motion group, the Bochner representation
 of these functions is an immediate consequence of  Harish-Chandra's theory (\cite{Hel}).
There are, however, no group-theoretical 
 interpretations known for the kernel $E_k$ so far.

\subsection{Heat kernels and  heat semigroups}

We start with a motivation:
Consider the following initial-value problem for the classical heat equation in $\b R^N$:
\begin{equation}
\begin{cases}\label{(4.500)} 
\Delta u - \partial_t u\, = 0 & \text{on $\b R^N\times(0,\infty)$},\\
     u(\,.\,,0) = f
   \end{cases}
\end{equation}
with initial data $f\in C_0(\b R^N)$, the space of continuous functions on $\b R^N$ which vanish
at infinity.  (We could equally take data 
from $C_b(\b R^N)$, but $C_0(\b R^N)$ is more convenient in the following considerations). 
The basic idea to solve \eqref{(4.500)}
is to carry out a Fourier transform with respect to $x$.  This yields the candidate
\begin{equation}\label{(4.501)} 
u(x,t) = g_t * f(x) \,=\, \int_{\b R^N} g_t(x-y)f(y) dy \quad (t>0),\end{equation}
where  $g_t$ is the Gaussian kernel 
\[ g_t(x) \,=\, \frac{1}{(4\pi t)^{n/2}}\, e^{-|x|^2/4t}. \]
It is a well-known fact from classical analysis that \eqref{(4.501)} is in fact the
unique bounded solution
 within the class $C^2(\b R^N\times (0,\infty))\cap C(\b R^N\times [0,\infty))$.

\begin{exercise}
Show that $H(t)f(x) := g_t*f(x)$ for $t>0, \,H(0):= id$ defines 
a strongly continuous contraction semigroup on the Banach space $(C_0(\b R^N), \|.\|_\infty)$ 
in the sense of the  definition given below. \\
\textit{Hint:} Once contractivity is shown, it suffices 
to check the continuity for functions from the
Schwartz space $\scr S(\b R^N).$ For this, use the Fourier inversion theorem.
\end{exercise}

\begin{definition}
Let $X$ be a Banach space. A one-parameter family $(T(t))_{t\geq 0}$ of bounded linear operators on $X$ is called a \textit{strongly continuous semigroup} on $X$, if it satisfies \parskip=-1pt
\begin{enumerate}
\item[\rm{(i)}]  $T(0) = id_X\,, \quad T(t+s) = T(t) T(s) \quad\text{for all }\,t,s\geq 0$
\item[\rm{(ii)}] The mapping $t\mapsto T(t)x$ is continuous on $[0,\infty)$ for all $x\in X$.
\end{enumerate}
A strongly continuous semigroup is called a \textit{contraction semigroup}, if 
$\|T(t)\|\leq 1$ for all $t\geq 0.$ 
\end{definition}

\noindent
Let $L(X)$ denote the space of bounded linear operators in $X$. 
If $A\in L(X),$ then 
\[ e^{tA} = \sum_{n=0}^\infty \frac{t^n}{n!} A^n\, \in L(X)\]
defines a strongly continuous semigroup on $X$ (this one is even continuous with respect to the
uniform topology on $L(X)$). We obviously have 
\[ A = \lim_{t\downarrow 0} \frac{1}{t} (e^{tA} -id) \quad\text{in }\, L(X).\]

\begin{definition} The \textit{generator} of a strongly continuous semigroup $(T(t))_{t\geq 0}$
in $X$ is defined by
\begin{align} Ax:=& \lim_{t\downarrow 0} \frac{1}{t} (T(t)x-x), \quad\text{with domain}\notag\\
   \mathcal D(A):=& \{x\in X:  \lim_{t\downarrow 0} \frac{1}{t} \bigl(T(t)x-x\bigr) 
\,\,\text{exists in }\, X
\}.\end{align}
\end{definition}

\begin{theorem} The generator $A$ of  $(T(t))_{t\geq 0}$ is densely defined and closed.
\end{theorem}

An important issue in the theory of operator semigroups and evolution equations 
are criteria which characterize
 generators of strongly continuous semigroups. 

Let us return  to the Dunkl setting. As before,  $\Delta_k$ denotes the Dunkl Laplacian
associated with a finite reflection group on $\b R^N$ and some multiplicity function $k\geq 0,$
and the index $\gamma$ is defined according to \eqref{(1.1a)}.
We are going to consider the following initial-value problem for the 
Dunkl-type heat operator $\Delta_k-\partial_t$:

Find $u\in C^2(\b R^N\times (0,\infty))$ which is continuous on 
$\b R^N\,\times [0,\infty)\,$ 
and satisfies
\begin{equation}
\begin{cases}\label{(4.21)} 
(\Delta_k - \partial_t)\, u\, =\, 0 & \,\text{on $\b R^N\times (0,\infty)$},\\
     u(\,.\,,0)\, =\, f & \text{$\in C_b(\b R^N).$}
   \end{cases}
\end{equation}

The solution of this problem is given, just as in the classical case $k=0$, in terms of a 
positivity-preserving semigroup.
We shall essentially follow the treatment of \cite{R2}.

\begin{lemma} The function 
\[ F_k(x,t) := \frac{1}{(2t)^{\,\gamma+N/2}c_k} e^{-|x|^2/4t}\]
solves the generalized  heat equation
 $\Delta_k u - \partial_t u = 0$ on $\b R^N\times 
(0,\infty).$ 
\end{lemma}

\begin{proof} A short calculation. Use the product rule \eqref{(1.2)} as well as the identity\\
$\sum_{i=1}^N T_i(x_i) = N+2\gamma\,.$
\end{proof}

\noindent
$F_k$  generalizes the fundamental solution for the classical heat equation
 which is given by $\, F_0(x,t) = g_t(x)$ (as defined above).
It is easily checked that
\[\int_{\b R^N} F_k(x,t)\, w_k(x) dx\,=\, 1 \quad \text{for all}\>\> t>0.\]

In order to solve \eqref{(4.21)}, it 
suggests itself to 
apply  the Dunkl transform under 
suitable decay assumptions on the initial data.
In the classical case, the heat kernel $g_t(x-y)$ on $\b R^N$ is obtained from 
the  fundamental solution simply by translations. In the Dunkl setting, it is 
still possible to define a generalized translation which matches the action of
the Dunkl transform, i.e. makes it a homomorphism on suitable function spaces. 

\noindent
The notion of a \textit{generalized 
translation}  in the Schwartz space $\mathscr S(\b R^N)$ is as follows (c.f. \cite{R2}):
\begin{equation}\label{(4.75)}
 \tau_y f(x):= \frac{1}{c_k} \int_{\b R^N} \widehat f^{\,k}(\xi) 
\,E_k(ix,\xi) E_k(iy, \xi)\, w_k(\xi) d\xi; \quad y\in \b R^N. 
\end{equation}
In the same way, this could be done in $L^2(\b R^N,w_k)$. A  powerful extension
to $C^\infty(\b R^N)$ is due to  Trim\`eche \cite{T}.
Note that in case $k=0$, we simply have $\,\tau_y f(x) = f(x+y).$ 
In the rank-one case, the above
 translation coincides with the convolution on a so-called signed
hypergroup structure which was defined in \cite{R1}; see also \cite{Ros}.  
Similar structures
are not yet known in higher rank cases.
   Clearly, $\tau_y f(x) = \tau_x f(y)$; 
moreover, the inversion theorem for the Dunkl transform assures that $\,\tau_0 f = f$ and 
\[ (\tau_y f)^{\wedge k}(\xi) 
= E_k(iy,\xi) \widehat f^{\,k}(\xi).\]
 From this it is not hard to see  that 
$\tau_y f$ belongs to $\mathscr S(\b R^N)$ again. 
Let us now consider the ``fundamental solution'' $F_k(.\,,t)$ for $t>0$. A 
short calculation, using the reproducing property Prop. \ref{P:Reprod}(2), shows that
\begin{equation}\label{(4.80)}
 \widehat F_k^k(\xi,t)\,=\, c_k^{-1}e^{-t|\xi|^2}.
\end{equation}
By  the quoted reproducing formula  one therefore obtains from 
\eqref{(4.75)} the representation
\begin{equation}
 \tau_{-y}F_k(x,t)\,=\, \frac{1}{(2t)^{\gamma+N/2}c_k} e^{-(|x|^2 + |y|^2)/4t}\,
E_k\Bigl(\frac{x}{\sqrt{2t}},\frac{y}{\sqrt{2t}}\Bigr).\notag
\end{equation}

\noindent
This motivates the following

\begin{definition}
The generalized heat kernel $\Gamma_k$ is defined by
\[ \Gamma_k(t,x,y):=\, \frac{1}{(2t)^{\gamma +N/2}c_k}\,e^{-(|x|^2 + |y|^2)/4t}\,
E_k\Bigl(\frac{x}{\sqrt{2t}},\frac{y}{\sqrt{2t}}\Bigr),\quad x,y\in \b R^N,\> 
t>0.\]
\end{definition}

\noindent
Notice in particular that $\Gamma_k >0$ (thanks to Corollary \ref{C:Bochner})
and  that $\,y\mapsto \Gamma_k(t,x,y)$ belongs to
$\mathscr S(\b R^N)$ for fixed $x$ and $t$.
We collect a series of further fundamental properties of this kernel which are all 
more or less straightforward.

\begin{lemma}\label{L:Heatprop} The heat kernel $\Gamma_k$ has the
  following 
properties: 
\begin{enumerate}\itemsep=1pt
\item[\rm{(1)}] $\,\displaystyle  \Gamma_k(t,x,y) \,=\, 
 c_k^{-2} \int_{\b R^N} e^{-t|\xi|^2}\, E_k(ix,\xi)\, 
E_k(-iy,\xi) \, w_k(\xi) d\xi\,.$ 
\item[\rm{(2)}] $\,\displaystyle \int_{\b R^N} \Gamma_k(t,x,y)\,w_k(y) 
dy\,=\,1.$
\item[\rm{(3)}] $\,\displaystyle \Gamma_k(t,x,y)\,\leq\, 
  \frac{1}{(2t)^{\gamma+N/2}c_k}\,\max_{g\in G} e^{-|gx-y|^2/4t}\,.$ 
\item[\rm{(4)}] $\,\displaystyle \Gamma_k(t+s,x,y)\,=\, \int_{\b R^N}
  \Gamma_k(t,x,z)\,\Gamma_k(s,y,z)\,w_k(z)dz.$
\item[\rm{(5)}] For fixed $y\in \b R^N$, the function $u(x,t):= 
\Gamma_k(t,x,y)$ solves the generalized heat equation $\, \Delta_k u =
\partial_t u\,$ on 
$ \b R^N\times(0,\infty)$.
\end{enumerate}
\end{lemma}

\begin{proof} (1) is clear from the definition of generalized translations.
 For details concerning (2) 
see \cite{R2}. (3) follows from
our estimates on $E_k$, while (4)
is obtained by inserting (1) for one of the kernels in the integral. Finally, 
(5) is obtained from differentiating (1) under the integral. For details see again \cite{R2}.
\end{proof}

\begin{definition}\label{D:Heat}
For $f\in C_b(\b R^N)$ and $t\geq 0$ set 
\begin{equation}\label{(3.16)}
 H(t)f(x):= \begin{cases}
    \displaystyle\int_{\b R^N} \Gamma_k(t,x,y) f(y)\,w_k(y)dy & 
    \text{if $\,\,t>0$},\\
     f(x) & \text{if $\,\,t=0$}
   \end{cases}
\end{equation}
\end{definition}

Notice that the decay of $\Gamma_k$ assures the convergence of the integral.
The properties  of the operators $H(t)$ are most easily described on the Schwartz space $\scr S(\b R^N)$. The following theorem is completely analogous to the classical case.

\begin{theorem} \label{T:Schwartz}
Let $f\in \scr S(\b R^N)$. Then $\, u(x,t):=  H(t)f(x)$ 
solves the initial-value problem \eqref{(4.21)}.
Moreover, $H(t)f$  has the following properties: 
  \parskip=-1pt
\begin{enumerate}
\item[\rm{(1)}]  $\displaystyle H(t)f\in \scr S(\b R^N)$ for all $t>0$.
\item[\rm{(2)}] $\displaystyle H(t+s)\,f \,=\, H(t)H(s)f\,$
  for all $s,t\geq 0$.
\item[\rm{(3)}] $\displaystyle \|H(t)f - f\|_{\infty,\b R^N} \, \to 0\> $ as
$t\to 0$. 
\end{enumerate}
\end{theorem}

\begin{proof} (Sketch)  By use of 
  Lemma \ref{L:Heatprop} (1)  and Fubini's theorem, we write
\begin{equation}\label{(4.12)}
 u(x,t)\,=\, H(t)f(x)\,=\, c_k^{-1} \int_{\b R^N} e^{-t|\xi|^2} 
\widehat f^{\,k}(\xi)\, E_k(ix,\xi)\,  w_k(\xi) d\xi \quad(t>0).
\end{equation}
 In view of the inversion theorem for the Dunkl transform, 
this holds for $t=0$ as well. Properties  (1) and (3) as well as the differential equation
 are now easy consequences.
Part (2) follows from the 
reproducing formula for $\Gamma_k$ (Lemma \ref{L:Heatprop} (4)).
\end{proof} 

\begin{exercise} Carry out the details in the proof of Theorem \ref{T:Schwartz}.
\end{exercise}

\noindent
We know that the heat kernel $\Gamma_k$ is positive; this implies that
$H(t)f \geq 0$ if $f\geq 0$. 

\begin{definition}
Let $\Omega$ be a locally compact Hausdorff space. 
A strongly continuous semigroup $(T(t))_{t\geq 0}$ on $(C_0(\Omega), \|.\|_\infty)$ 
is called a \textit{Feller-Markov semigroup}, if 
it is contractive and positive, i.e. $\, f\geq 0$ on $\Omega$ implies that $T(t)f \geq 0$ on 
$\Omega$ for all $t\geq 0.$
\end{definition}

We shall prove that
 the linear operators $H(t)$ on $\scr S(\b R^N)$ extend to a Feller-Markov semigroup
on the Banach space $(C_0(\b R^N),\|.\|_\infty)$. 
This could be done by direct calculations similar to the usual procedure for the
classical heat semigroup, relying on the positivity of the kernel $\Gamma_k$. 
 We do however prefer to give a proof which does not require this rather deep result,
but works on the level of the tentative generator. The tool is the 
following useful variant of the Lumer-Phillips theorem, which  characterizes
 Feller-Markov semigroups 
in terms of a ``positive maximum principle'', see e.g. \cite{Kal}, Thm. 17.11.
In fact, this Theorem motivated the positive minimum principle
\ref{T:Posmin}  in the positivity-proof for $V_k$.

\begin{theorem}\label{T:Feller}
Let $(A, \mathcal D(A))$ be a 
densely defined linear operator in $(C_0(\Omega), \|.\|_\infty)$. Then $A$ is closable,
and its closure $\overline A$ generates a Feller-Markov semigroup on $C_0(\Omega),$ if and only
if the following conditions are satisfied:
\begin{enumerate}
\item[\rm{(i)}] If $f\in \mathcal D(A)$ then also $\overline f\in \mathcal D(A)$ and $A(\overline f) 
= \overline{A(f)}.$
\item[\rm{(ii)}] The range of $\lambda id -A$ is dense in $C_0(\Omega)$ for some $\lambda >0$.
\item[\rm{(iii)}] If $f\in \mathcal D(A)$ is real-valued with a non-negative maximum in $x_0\in\Omega$, i.e. \\ $0\leq f(x_0) = \max_{x\in\Omega} f(x),$ then $Af(x_0) \leq 0.$ (Positive 
maximum principle).
\end{enumerate}
\end{theorem}

We consider the Dunkl Laplacian $\Delta_k$ as a densely 
defined linear operator in $C_0(\b R^N)$ with domain $\mathscr S(\b R^N)$. The following Lemma implies that it satisfies the positive maximum principle:

\begin{lemma}
Let $\Omega\subseteq\b R^N$ be open and $G$-invariant. 
If a real-valued function $f\in C^2(\Omega)$ attains an absolute maximum at  
$x_0\in\Omega$, i.e. $f(x_0)= \sup_{x\in\Omega} f(x)$, then
\[\Delta_k f(x_0)\,\leq\, 0\,.\]
\end{lemma}

\begin{exercise} \label{E:mini}
Prove this lemma in the case that $\langle\alpha ,x_0\rangle \not=0$ for all 
$\alpha\in R$. (If $\langle\alpha ,x_0\rangle =0$ for some $\alpha\in R$, one has to argue 
more carefully; for details see \cite{R2}.)
\end{exercise}

\begin{theorem}\label{T:Fellermain}
The operators $(H(t))_{t\geq 0}$ define a Feller-Markov semigroup on $C_0(\b R^N)$. 
Its generator is the closure $\overline \Delta_k$ of $(\Delta_k, \mathscr S(\b R^N))$. 
This semigroup is called the \textit{generalized heat semigroup} on $C_0(\b R^N)$.
\end{theorem}

\begin{proof} In the first step, we check that $\Delta_k$ (with domain $\mathscr S(\b R^N)$)
satisfies the conditions of Theorem \ref{T:Feller}: Condition (i) is obvious and (iii) is an immediate consequence of the previous lemma. 
Condition (ii) is 
also satisfied, because $\lambda id  -\Delta_k$ maps $\mathscr S(\b R^N)$ 
onto itself for each $\lambda >0$;  this follows from the fact that the 
Dunkl transform  is a  homeomorphism of $\mathscr S(\b R^N)$ and 
$\,\bigl((\lambda I  -\Delta_k)f\bigr)^{\wedge k}(\xi) = 
(\lambda +|\xi|^2)\widehat f^{\,k}(\xi)$. 
Theorem \ref{T:Feller} now implies that $\Delta_k$ is closable, and that its closure $\overline\Delta_k$ generates a Feller-Markov semigroup $(T(t))_{t\geq 0}$. 
It remains to show that $T(t) = H(t)$ on $C_0(\b R^N)$. Let first $f\in \scr S(\b R^N).$ 
From basic facts in semigroup theory, it follows that  the function 
$t\mapsto T(t)f$ is the unique solution of the so-called abstract Cauchy problem
\begin{equation}
\begin{cases}\label{(4.70)} 
\displaystyle\frac{d}{dt}u(t)= \overline\Delta_k u(t) & \text{for $t>0$},\\
     u(0) = f & \text{}
   \end{cases}
\end{equation}
within the class of all (strongly) continuously differentiable functions $u$ 
on $[0,\infty)$ with values in  $(C_0(\b R^N), \|.\|_\infty).$ It is easily seen from Theorem 
\ref{T:Schwartz}, and in particular from formula \eqref{(4.12)}, that
$t\mapsto H(t)f$ satisfies these conditions. Hence $T(t)= H(t)$ on $\scr S(\b R^N)$.
This easily implies that $\Gamma_k \geq 0$ (which we did not presuppose for the proof!),
 and therefore the operators $H(t)$ are also contractive
on $C_0(\b R^N)$.  A density argument now finishes the proof. 
\end{proof}

Based on this result, it is checked by standard arguments that for data $f\in C_b(\b R^N),$
the function $u(x,t):= H(t)f(x)$ solves the initial-value problem \eqref{(4.21)}. Uniqueness
results are established by means of maximum principles, just as with the classical heat equation.
Moreover, 
the heat semigroup $(H(t))_{t\geq 0}$  can also be defined (by means of \eqref{(3.16)})
on the Banach spaces $L^p(\b R^N, w_k), 1\leq p <\infty$. In case $p=2$, the following is
easily seen by use of the Dunkl transform: 

\begin{proposition}\label{L2} \cite{R2}
The operator $(\Delta_k, \mathscr S(\b R^N))$ in $L^2(\b R^N, w_k)$ 
is densely defined and closable. Its closure generates a strongly continuous 
and positivity-preserving  contraction
semigroup on $L^2(\b R^N, w_k)$ which is given by
\[ H(t)f(x)= 
  \,\int_{\b R^N} \Gamma_k(t,x,y) f(y)w_k(y)dy\,,\quad (t>0).\]
\end{proposition}

Theorem \ref{T:Fellermain}
 was the starting point in \cite{RV3} to construct an associated Feller-Markov  process on $\b R^N$  which can be
 considered a generalization of the usual Brownian motion. The transition probabilities 
of this process are defined in terms of a semigroup of Markov kernels of $\b R^N$,  as follows:
For $x\in \b R^N$ and a Borel set $\,A\in \mathscr B(\b R^N)$ put
\[ P_t(x,A):=\, \int_A \Gamma_k(t,x,y)w_k(y)dy \quad (t>0), \quad
P_0(x,A) := \delta_x(A),\]
with $\delta_x$ denoting the point measure in $x\in \b R^N$. Then
$(P_t)_{t\geq 0} \,$ is a semigroup of Markov kernels on $\b R^N$ in
the following sense:\parskip=-1pt
\begin{enumerate}
\item[\rm{(1)}] Each $P_t$ is a Markov kernel, and for all $s,t\geq 0,
  \> x\in \b R^N$ and $A\in \mathscr B(\b R^N)$,
\[ P_s\circ P_t (x,A):=\, \int_{\b R^N} P_t(z,A)\,P_s(x,dz)\,\,=\,
P_{s+t}(x,A).\]
\item[\rm{(2)}] The mapping $[0,\infty)\,\to M^1(\b R^N), \,
  t\mapsto P_t(0,\,.\,),$ is continuous with respect to the 
$\sigma(M^1(\b R^N, C_b(\b R^N))$-topology.
\end{enumerate}
\noindent
Moreover, the semigroup $(P_t)_{t\geq 0} \,$ has the following
particular property:
\begin{enumerate} 
\item[\rm{(3)}] $\displaystyle   P_t(x,.\,)^{\wedge k}(\xi) =
  E_k(-ix,\xi)\, P_t(0,\,.\,)^{\wedge k}(\xi)\,$ for all $ \xi
  \in \b R^N,$
\end{enumerate}
hereby the Dunkl transform of the probability measures
$\,P_t(x,\,.\,)$ is defined by
\[ P_t(x,\,.\,)^{\wedge k}(\xi):= \int_{\b R^N} E_k(-i\xi,x)\,P_t(x,d\xi).\]
The proof of  (1) -- (3) is straightforward by the properties of
$\Gamma_k$ and Theorem \ref{T:Fellermain}.

\noindent
In the classical case $k=0$, property (3) is equivalent to $(P_t)_{t\geq 0}$ being 
translation-invariant, i.e. 
\[P_t(x+y, A+y)\,=\, P_t(x,A) \quad\text{for all }y\in \b R^N.\] 
In our general setting, a positivity-preserving translation on $M^1(\b R^N)$ cannot be expected
(and does definitely not exist in the rank-one case according to \cite{R1}). Property (3) 
thus serves as a substitute for translation-invariance.  The reader can see
\cite{RV3} for a 
study of the semigroup $(P_t)_{t\geq 0} \,$ and  the associated Feller-Markov process.

\section{Asymptotic analysis for the Dunkl kernel}

This final section deals with the asymptotic behavior of the Dunkl kernel $E_k$ with $k\geq 0$ 
 when one of its
arguments is fixed and the other tends to infinity either
 within a Weyl chamber 
of the associated reflection group, or within a suitable complex domain. These results
are contained in \cite{dJR}.
They
 generalize the well-known asymptotics of the confluent hypergeometric
function $\phantom{}_1F_1$  to
the higher-rank setting. One motivation to study the asymptotics of $E_k$ 
is to determine the asymptotic behavior of the Dunkl-type heat kernel $\Gamma_k$ for short times. 
Partial results in this direction  were obtained in \cite{R4}.

Recall from Prop.~\ref{L2} that $\Gamma_k$ is the kernel of the generalized heat semigroup in the weighted
 $L^2(\b R^N, w_k)$. We want to compare it with the free Gaussian kernel $\Gamma_0$. For this,
it is appropriate to transfer the semigroup $(H(t))_{t\geq 0}$ from  $L^2(\b R^N, w_k)$
to the unweighted  space $L^2(\b R^N)$, which leads to the strongly
continuous contraction semigroup
\[ \widetilde H(t)f:= w_k^{1/2} H(t)\bigl(w_k^{-1/2}f), \quad f\in
L^2(\b R^N).\]
The corresponding renormalized heat kernel is given by 
\[ \widetilde \Gamma_k(t,x,y):=\, \sqrt{ w_k(x)w_k(y)}\,
\Gamma_k(t,x,y).\]
The generator of $(\widetilde H(t))_{t\geq 0}$ is 
 the  gauge-transformed version of the Dunkl Laplacian discussed
in connection with CMS-models, 
\[
\mathcal F_k\,=\, \Delta -2\sum_{\alpha\in R_+}
 \frac{k(\alpha)}{\langle\alpha,x\rangle^2} (k(\alpha)
 -\sigma_\alpha) 
\] 
(with suitable domain).
$\mathcal F_k$ can be considered  a perturbation of the 
Laplacian $\Delta$. This suggests  that 
                   within the Weyl chambers 
of $G$, the heat kernel
$\widetilde\Gamma_k(t,x,y)$ should not ``feel'' the reflecting hyperplanes and 
behave for short times like the free
Gaussian  kernel 
$\,  \Gamma_0(t,x,y)\,=\, g_t(x-y)$, in other words, we have the conjecture
\begin{equation}\label{(7.1)}
 \lim_{t\downarrow 0} \frac{\sqrt{w_k(x)w_k(y)}\,\Gamma_k(t,x,y)}{\Gamma_0(t,x,y)}\,=\,1
\end{equation}
provided $x$ and $y$ belong to the same (open) Weyl chamber.
In \cite{R4}, this could be proven true  only for a restricted
range of arguments $x,y,$ and by rather technical methods (completely different from those below).

\begin{example} \textbf{The rank-one case.}
Here $E_k$ is explicitly known. According to  Ex. \ref{E:1kern}, 
\[E_k (z,w)\,=\, e^{zw}\cdot
\phantom{}_1F_1(k,2k+1,-2zw), \quad z,w\in \b C.\]
 The  confluent hypergeometric function $\phantom{}_1F_1$ has well-known asymptotic expansions in the sectors 
\[ S_+ = \{z\in \b C: -\pi/2 < \text{arg}(z) < 3\pi/2\}, \quad
   S_- = \{z\in \b C: -3\pi/2 < \text{arg}(z) < \pi/2\},\]
see for instance \cite{AS}. They are of the form
\[\phantom{}_1F_1(k,2k+1,z) \,=\, \frac{\Gamma(2k+1)}{\Gamma(k)} e^z z^{-k-1}
 \bigl(1 + \mathcal O(\frac{1}{|z|})\bigr) \,+\, \frac{\Gamma(2k+1)}{\Gamma(k+1)} 
   e^{\pm i\pi k} z^{-k} \bigl(1+ \mathcal O(\frac{1}{|z|})\bigr), \]
with $\pm$ for $z\in S_{\pm}$. Specializing to the right half plane
\[ H = \{z\in \b C: \text{Re} z\geq 0\}\] 
we thus obtain
\[\lim_{zw\to\infty, zw\in H} (zw)^k e^{-zw} E_k(z,w) \,=\, \frac{\Gamma(2k+1)}{2^k\Gamma(k+1)}.\]
\end{example}

Let us now turn to the general case of an arbitrary reflection group and 
multiplicity parameter $k\geq 0$. 
Let $C$
denote the Weyl chamber attached with the positive subsystem $R_+$,
\[C= \{x\in \b R^N: \langle
\alpha,x\rangle >0 \,\,\text{ for all }\, \alpha\in R_+\},\]
and for $\delta >0,$
\[C_\delta :=
\{x\in C: \langle\alpha,x\rangle\,>\,
\delta |x| \,
  \text{ for all }\,\alpha\in R_+\}.\]
The main result given here is the following  asymptotic
behavior, uniform for the variable tending to infinity in cones $C_\delta$:

\begin{theorem}\label{T:main1}
There exists a constant non-zero vector $v=(v_g)_{g\in G}\in\b C^{|G|}$ such
that for all $y\in C, \,g\in G$ and each $\delta>0,$
\[\lim_{|x|\to\infty,\,x\in C_\delta} {\sqrt{w_k(x)w_k(y)}}\,e^{-i\,\langle x,gy\rangle} E_k(ix,gy)\,=\, v_g.\]
\end{theorem}

\noindent
Notice that one variable is being fixed. A locally uniform result with respect to both 
variables should be true, but is open yet. Also, the explicit values of the constants $v_g$ --
apart from $v_e$ -- are not known. We come back to this point later.
An immediate consequence of Theorem \ref{T:main1} is the
following ray asymptotic for the Dunkl kernel (already conjectured in  \cite{Du3}):

\begin{corollary}\label{C:ray}
For all $x,y\in C$ and $g\in G$,
\[ \lim_{t\to\infty} t^\gamma \,e^{-it\,\langle x,gy\rangle} E_k(itx,gy)\,=\,
   \frac{v_g}{\sqrt{w_k(x)w_k(y)}}\,,\]
the convergence being locally uniform with respect to the parameter $x$.
\end{corollary}

In the particular case $g=e$ (the unit of $G$), this latter result
can be extended to a larger range of
complex arguments by use of the Phragm\' en-Lindel\"of principle for the
right half plane $H$ (see e.g. \cite{Ti}):

\begin{proposition}\label{P:Phragmen}
 Suppose $f: H\to \b C$ is  analytic and regular in $H\cap\{z\in \b C: |z|> R\}$ for some $R>0$  with
$\lim_{t\to\infty}f(it) = a,\,\lim_{t\to\infty}f(-it) =b$ and  for each $\delta >0,$ 
\[ f(z) = \mathcal O\bigl(e^{\delta|z|}\bigr) \quad\text{as }\, z\to\infty 
\,\,\text{within }H.\]
Then $a=b$ and $f(z)\to a$ uniformly as $z\to\infty$ in $H$. 
\end{proposition}

\begin{theorem}\label{T:main3} Let $x,y\in C$. Then
\[ \lim_{z\to\infty, z\in H}
z^\gamma e^{-z\langle x,y\rangle} E_k(zx,y)\,=\,\frac{i^\gamma v_e}{\sqrt{w_k(x)w_k(y)}}.
\]
Here  $z^\gamma$ is the holomorphic branch in
$\b C\setminus
\{x\in \b R: x\leq 0\}$ with $1^\gamma = 1$.
\end{theorem}

\begin{proof}  Consider
\[ G(z):= z^\gamma e^{-z\langle x,y\rangle}E_k(zx,y).\]
The estimate of Prop. \ref{P:Roughestim} on $E_k$ implies that
$G$ satisfies the required growth 
bound in Prop. \ref{P:Phragmen} (here it is of importance that $x$ and $y$ lie in the
 same Weyl chamber), and Corollary \ref{C:ray}
 assures that $G$ has limits along the boundary
lines of $H$. 
\end{proof}

When restricted to real arguments, Theorem \ref{T:main3} implies the above
stated short-time asymptotic for the heat kernel $\Gamma_k$:

\begin{corollary}\label{C:heat} For all $x,y\in C$,
\[\lim_{t\downarrow 0} \frac{\sqrt{w_k(x)\,w_k(y)}\,\,
\Gamma_k(t,x,y)}{\Gamma_0(t,x,y)}\,=\,1\,.\]
\end{corollary}

Hereby the precise value of the limit
follows 
from the results of \cite{R4}. Along with it, we thus obtain the value of $v_e$: 
\[ v_e \,=\, i^{-\gamma}\frac{c_k}{c_0}.\]
We shall now give an outline of the proof of   Theorem \ref{T:main1}. It is
based on the analysis of an associated  system of
first order  differential equations, which is derived
from the eigenfunction characterization
 \eqref{(1.106)} 
 of $E_k$.  This approach goes
 back to \cite{dJ1}, where it was used to obtain
exponential estimates for the Dunkl kernel. 
Put 
\[\b R^N_{reg}:= \b R^N\setminus\{\langle\alpha\rangle^\perp, \alpha\in R\}\]
and  define
\[ \phi(x,y) = \,\sqrt{w_k(x)w_k(y)}\,
e^{-i\langle x,y\rangle} E_k(ix,y), \quad x,y\in \b R^N.\]
Observe that $\phi$  is symmetric in its arguments. We have to study
the asymptotic behavior of $x\mapsto \phi(x,y)$ along curves in $C,$
with the second component $y\in \b R^N_{reg}$ being
fixed. 
Let us introduce the auxiliary
 vector field
  $\,F = (F_g)_{g\in G}$  on $\b R^N\times\b R^N$ by
\[ F_g(x,y)\,:=\, \phi(x,gy).\]
For fixed $y$, we consider $F$ along a differentiable curve $\kappa: (0,\infty)\to C$.
The eigenfunction characterization 
of $E_k$ then  translates into a first order ordinary differential
equation for $t\mapsto F(\kappa(t),y)$. 
Below, we shall determine the asymptotic behavior of its solutions,
provided $\kappa$ is admissible in the following sense:

\begin{definition} A $C^1$-curve\\
\begin{minipage}[t]{.51\linewidth}
$\kappa: (0,\infty)\to C$ is called \emph{admissible}, if
it satisfies the subsequent conditions:
\parskip=-1pt
\begin{enumerate}
\item[\rm{(1)}] There exists a constant $\delta >0$ such that
   $\kappa(t)\in C_\delta$ for all $t>0.$
\item[\rm{(2)}] $\lim_{t\to\infty} |\kappa(t)| \,=\infty\,$
 and $\kappa^\prime(t)\in C$ for all $t>0$.
\end{enumerate}
\end{minipage}\hfill
\begin{minipage}[t]{.49\linewidth}
\vspace*{-34pt}
\psfull
$$\epsfbox{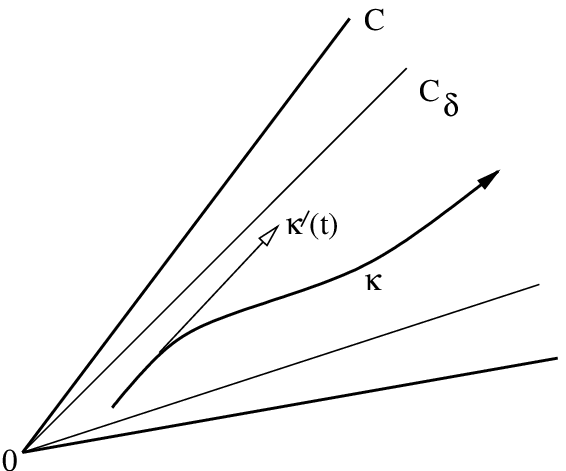}$$
\end{minipage}\hfill
\end{definition}

An important class of admissible curves are the rays $\kappa(t) = tx$ with
some fixed $x\in C$.  In a first step, it is shown that
 $t\mapsto F(\kappa(t),y)$ is asymptotically constant as
$t\to\infty$ for arbitrary admissible curves:

\begin{theorem} \label{T:Main22}  If  $\kappa: (0,\infty)\to C$ is
admissible, then for
every   $y\in C,$  the limit
\[\lim_{t\to\infty} F(\kappa(t),y)\]
exists in $\b C^{|G|}$, and is different from $0$.
\end{theorem}

\begin{proof} (Sketch)
 It is easily calculated that \eqref{(1.106)} translates into the
following
system of differential equations for $(F_g)_{g\in G}$, where
  $\xi, \,y\in \b R^N$ are fixed:
\begin{equation}\label{(7.2)} \partial_{\xi} F_g(x,y)\,=\, \sum_{\alpha\in R_+} k(\alpha)\,
  \frac{\langle\alpha,\xi\rangle}{\langle\alpha,x\rangle}\,
e^{-i\langle \alpha,\,x\rangle\langle\alpha,\,gy\rangle}\cdot
   F_{\sigma_\alpha g}(x,y)\quad (x\in \b R^N_{reg}).
\end{equation}

From this, one obtains a differential equation for $x(t):=  F(\kappa(t),y)$ of the form 
$x^\prime(t) = A(t)x(t)$, with a continuous matrix function 
 $A: (0,\infty)\to \b C^{|G|\times |G|}\,$. The proof of Theorem \ref{T:Main22} is accomplished by
 verifying that $A$ 
satisfies the conditions of the following classical theorem on the asymptotic integration of ordinary differential equations (hereby of course, the admissibility conditions
on $\kappa$ come in).
\end{proof}

\begin{theorem}\label{P:Wintner} \cite{East},\,\cite{W}.
Consider the
linear differential equation
\begin{equation}\label{(1.1)}
 x^\prime(t)\,=\, A(t) x(t),
\end{equation}
where  $A:[t_0,\infty)\to \b C^{n\times n}$ is a continuous matrix-valued
function satisfying
the following  integrability conditions:
\parskip=-1pt
\begin{enumerate}
\item[\rm{(1)}] The matrix-valued improper Riemann
   integral $\,\int_{t_0}^\infty A(t)dt\,$ converges.
\item[\rm{(2)}] $ t\mapsto A(t)\int_{t}^\infty A(s)ds\,$ belongs to  
$L^1\bigl([t_0,\infty), \b C^{\,n\times n}\bigr).$
\end{enumerate}
Then $\eqref{(1.1)}$ has a fundamental system $\Phi$ of solutions 
which  satisfies$\, \lim_{t\to\infty}\Phi(t) \,=\, Id.$
\end{theorem}

Notice that in the situation of this theorem, for each solution $x$ 
of $x^\prime(t) = A(t) x(t)$ the limit $\lim_{t\to\infty} x(t)$ exists, and
is different from zero, unless $x\equiv 0$.

\begin{proof}[Proof of Theorem \ref{T:main1}]
It remains to
  show that
the limit value according to Theorem \ref{T:Main22}
is actually independent of $y$ and $\kappa$; the assertion
 is then  easily obtained. The stated independence is 
accomplished  as follows: 
In a first step, we show that there exists a non-zero  vector
  $v(y) = (v_g(y))_{g\in G}\in \b C^{|G|}$ such that for each admissible  $\kappa$,
\begin{equation}\label{(8.3)} \lim_{t\to\infty} F(\kappa(t),y)\,=\,v(y).
\end{equation}
This can be achieved via  interpolation of the 
 admissible curves $\kappa_1\,,\,\kappa_2$ by  a third  admissible  curve $\kappa$,
such that  equality of all three limits is being enforced.
Next, we focus on admissible rays. 
Observe that $\,F_g(tx,y) = F_{g^{-1}}(ty,x)$ for all $g\in G$ and
$x,y\in C$. Together with
\eqref{(8.3)}, this implies that
 $\, v_g(y) = v_{g^{-1}}(x)$, and therefore
also
 $\,v_g(x)\, = \,v_{g^{-1}}(x)\,=\, v_g(y)\, =: v_g.$  Put $v= (v_g)_{g\in G}$. Then
\[ \lim_{t\to\infty} F(\kappa(t),y) = v\]
for every admissible $\kappa$ and every $y\in C$.
\end{proof}

The asymptotic result of Theorem \ref{T:main1} also allows
 to deduce at least a certain amount of information about
 the structure of the intertwining operator $V_k$ and its representing measures
 $\mu_x^k$  according to formula \eqref{(4.11)}. The key for
our approach
 is the following
simple observation: according to Corollary \ref{C:Bochner}, one may write
\[
 E_k(x,-i\xi) =\, \int_{\b R^N}
e^{-i\langle \xi, y\rangle}d\mu_x^k(y)\,=\, \widehat{\mu_x^k}(\xi) \quad
(x,\xi\in \b R^N),
\]
where $\widehat \mu$ stands for the (classical)
Fourier-Stieltjes transform of  $\mu\in M^1(\b R^N)$,
\[ \widehat\mu(\xi) \,=\, \int_{\b R^N} e^{-i\langle \xi,y\rangle} d\mu(y).\]
Recall that a measure $\mu\in M^1(\b R^N)$ is called continuous, if $\mu(\{x\}) = 0$
for all $x\in\b R^N$.  There is a  well-known criterion of Wiener   
which characterizes Fourier-Stieltjes transforms
of continuous measures on
locally compact abelian groups, here $(\b R^N, +);$
 see for instance Lemma 8.3.7 of \cite{GMG}:

\begin{lemma} (Wiener) \label{L:Wiener}
For  $\mu\in M^1(\b R^N)$ the following
properties are equivalent: \parskip=-2pt
\begin{enumerate}\itemsep=-1pt
\item[\rm{(1)}] $\mu$ is continuous.
\item[\rm{(2)}] $\displaystyle
\lim_{n\to\infty} \frac{1}{n^N} \int_{\{\xi\in \b R^N:\,|\xi|\leq n\}} |\widehat \mu(\xi)|^2
d\xi\,=\, 0.$
\end{enumerate}
\end{lemma}

\noindent
This yields the following result:

\begin{theorem}\label{T:main6} Let $k\geq 0$. Then apart from the case $k=0$ (i.e.  the
classical Fourier case),
the measure $\mu_x^k$ is continuous
for all $x\in \b R^N_{reg}$.
\end{theorem}

\noindent
We conclude with two open problems:
\begin{enumerate}
\item[\rm{(a)}]
  In the situation of the last theorem, prove that the measures
$\mu_x^k$ are even absolutely
continuous with respect to Lebesgue measure, provided $\{\alpha\in R: k(\alpha)> 0\}$ spans $\b R^N$. 
\item[\rm{(b)}] Determine the values of the constants $v_g\,,\, g\in G$.
\end{enumerate}

\section{Notation}

\noindent
We denote by $\b Z,\, \b R$ and $\b C$ the sets of integer, real and complex numbers respectively. Further, $\, \b Z_+ = \{n\in \b Z: n\geq 0\}$.
For a locally compact Hausdorff space $X$, we denote by
 $\,C(X), C_b(X), C_c(X), \, C_0(X)$  the spaces of continuous complex-valued functions on $X,$ those which are bounded, those with compact support, and those which vanish at infinity, respectively. 
Further, $M_b(X), \, M_b^+(X),\, M^1(X)$ are the spaces of regular bounded Borel measures
on $X,$ those which are positive, and those which are probability-measures, respectively. Finally, 
$\scr B(X)$ stands for the $\sigma$-algebra of Borel sets on $X$.


\begin{thebibliography}{9999}


\bibitem[AS]{AS}
Abramowitz, M.,  Stegun, I.A., 
\emph{Pocketbook of Mathematical
    Functions.} Verlag Harri Deutsch, Frankfurt/Main, 1984.

\bibitem[BF1]{BF1} Baker, T.H.,  Forrester, P.J., The Calogero-Sutherland 
model and generalized classical polynomials. \emph{Comm. Math. Phys.} 
188 (1997), 175--216. 
\bibitem[BF2]{BF2} Baker, T.H., Forrester, P.J.,  The Calogero-Sutherland 
model and polynomials with prescribed symmetry. \emph{Nucl. Phys. B}
492 (1997), 682--716. 
\bibitem[BF3]{BF3} Baker, T.H.,  Forrester, P.J., Non-symmetric Jack 
polynomials and integral kernels. \emph{Duke Math. J.} 95 (1998), 1--50.
\bibitem[BHKV]{BHKV} Brink, L.,  Hansson, T.H., Konstein, S., 
Vasiliev, M.A., The Calogero model - anyonic representation, fermionic
extension and supersymmetry. \emph{Nucl. Phys. B} 401 (1993), 591--612.
\bibitem[Ca]{Ca} Calogero, F., Solution of the one-dimensional N-body problems 
   with quadratic and/or inversely quadratic pair potentials.
\emph{J. Math. Phys.} 12 (1971), 419--436. 
\bibitem[Chi]{Chi} Chihara, T.S.,  \emph{An Introduction to Orthogonal 
Polynomials.} Gordon and Breach, 1978.
\bibitem[Da]{Da} Davies, E.B.: \emph{Spectral Theory and Differential
  Operators.} Cambridge University Press, 1995.
\bibitem[vD]{vD} van Diejen, J.F., Confluent hypergeometric orthogonal 
  polynomials related to the rational quantum Calogero system with harmonic 
  confinement. \emph{Comm. Math. Phys.} 188 (1997), 467--497.
\bibitem[DV]{DV} van Diejen, J.F., Vinet, L.
\textit{Calogero-Sutherland-Moser Models.}
CRM Series in Mathematical Physics, Springer-Verlag,  2000.
\bibitem[D1]{Du0}  Dunkl, C.F., Reflection groups and orthogonal 
polynomials on the sphere. \emph{Math. Z.}  197 (1988), 33--60.
\bibitem[D2]{Du1} Dunkl, C.F., Differential-difference operators
associated to reflection groups. \emph{Trans. Amer. Math. Soc.} 311 (1989),
 167--183.
\bibitem[D3]{Du2a} Dunkl, C.F., Operators commuting with Coxeter group
  actions on polynomials. In: Stanton, D. (ed.), \emph{Invariant Theory and
  Tableaux}, Springer,  1990, pp. 107--117.
\bibitem[D4]{Du2} Dunkl, C.F., Integral kernels with reflection group 
invariance.
 \emph{Canad. J. Math.} 43 (1991), 1213--1227.
\bibitem[D5]{Du3} Dunkl, C.F.,  Hankel transforms associated to finite 
reflection groups. In: \emph{Proc.
of the special session on hypergeometric functions on domains of 
positivity, Jack polynomials and applications.} Proceedings, Tampa 1991, 
Contemp. Math.  138 (1992), pp. 123--138. 
\bibitem[D6]{Du4} Dunkl, C.F., Intertwining operators associated to
  the group $S_3$. \emph{Trans. Amer. Math. Soc.} 347 (1995), 3347--3374.
\bibitem[D7]{Du7} Dunkl, C.F., Symmetric Functions and $B_N$-invariant spherical harmonics. Preprint; math.CA/0207122.
\bibitem[DJO]{DJO} Dunkl, C.F., de Jeu  M.F.E., Opdam, E.M., Singular 
polynomials for finite reflection groups. \emph{Trans. Amer. Math. Soc.} 346
(1994), 237--256.
\bibitem[DX]{DX} Dunkl, C.F., Xu, Yuan, \emph{Orthogonal Polynomials of Several Variables;}
    Cambridge Univ. Press, 2001.
\bibitem[E]{East} Eastham, M.S.P., \emph{The Asymptotic Solution of Linear
 Differential Systems: Applications of the Levinson Theorem.} Clarendon Press,
Oxford 1989.
\bibitem[FD]{FD} Fell, J.M.G., Doran, R.S., \emph{Representations of
  $*$-Algebras, Locally Compact Groups , and Banach-$*$-Algebraic
  Bundles}, Vol. 1. Academic Press,  1988.
\bibitem[GM]{GMG} C. Graham and O.C. McGehee, \emph{Essays in Commutative
Harmonic Analysis,} Springer Grundlehren 238, Springer-Verlag, New York 1979.
\bibitem[GB]{GB} Grove, L.C., Benson, C.T., \emph{Finite Reflection
    Groups;} Second edition. Springer, 1985.
\bibitem[Ha]{Ha} Ha, Z.N.C., Exact dynamical correlation functions of
  the Calogero-Sutherland model and one dimensional fractional
  statistics in one dimension: View from an exactly solvable
  model. \emph{Nucl. Phys. B} 435 (1995), 604--636. 
\bibitem[Hal]{Hal} Haldane, D., Physics of the ideal fermion gas:
  Spinons and quantum symmetries of the integrable Haldane-Shastry
  spin chain. In: A. Okiji, N. Kamakani (eds.), 
\emph{Correlation effects in low-dimensional electron
  systems.} Springer, 1995, pp. 3--20. 
\bibitem[He1]{He2}  Heckman, G.J.,  A remark on the Dunkl 
differential-difference operators. In: Barker, W.,  Sally, P. (eds.) 
\emph{Harmonic analysis on reductive groups.} Progress in Math. {\bf 101},
    Birkh\"auser, 1991. pp. 181 -- 191. 
\bibitem[He2]{He}  Heckman, G.J., Dunkl operators. S\'{e}minaire Bourbaki 828, 1996--97; 
\emph{Ast\'{e}risque} 245 (1997), 223--246. 
\bibitem[Hel]{Hel} Helgason, S.,  \emph{Groups and Geometric Analysis.} American
Mathematical Society, 1984.
\bibitem[Hu]{Hu} Humphreys, J.E., \emph{Reflection Groups and Coxeter
  Groups.} Cambridge University Press, 1990.
\bibitem[dJ1]{dJ1} de Jeu, M.F.E.,  The Dunkl transform. \emph{Invent. Math.} 
113 (1993), 147--162.
\bibitem[dJ2]{dJ3} de Jeu, M.F.E., Dunkl operators. Thesis, University
    of Leiden, 1994.
\bibitem[dJ3]{dJ4} de Jeu, M.F.E., Subspaces with equal closure. Preprint. math.CA/0111015.
\bibitem[K]{K} Kakei, S., Common algebraic structure for the 
 Calogero-Sutherland models. \emph{J. Phys. A} 29 (1996), L619--L624.
\bibitem[Kal]{Kal} Kallenberg, O., \emph{Foundations of Modern Probability.} Springer-Verlag, 
1997.
\bibitem[KS]{KS} Knop, F.,  Sahi, S.,  A recursion and combinatorial 
formula for Jack polynomials. \emph{Invent. Math.} 128 (1997), 9--22.   
\bibitem[LV]{LV} Lapointe  L., Vinet, L., Exact operator solution of the 
Calogero-Sutherland model. \emph{Comm. Math. Phys.} 178 (1996),
425--452.
\bibitem[La1]{L1} Lassalle, M., Polyn\^{o}mes de Laguerre 
g\'{e}n\'{e}ralis\'{e}s. \emph{C.R. Acad. Sci. Paris} t. 312 
S\'{e}rie I (1991), 725--728
\bibitem[La2]{L2} Lassalle, M., Polyn\^{o}mes de Hermite 
g\'{e}n\'{e}ralis\'{e}s. \emph{C.R. Acad. Sci. Paris} t. 313 
S\'{e}rie I (1991), 579--582.
\bibitem[M1]{M1} Macdonald, I.G., Some conjectures for root systems. \emph{SIAM J. Math. Anal.}
13 (1982), 988--1007.
\bibitem[M2]{M2} Macdonald, I.G., The Volume of a Compact Lie Group. 
\emph{Invent. Math.} 56 (1980), 93--95.
\bibitem[Mo]{Mo} Moser, J., Three integrable Hamiltonian systems connected with
 isospectral deformations, \emph{Adv. in Math.} 16 (1975), 197--220.
\bibitem[OP1]{OP1} Olshanetsky, M.A., Perelomov, A.M., Completely
  integrable Hamiltonian systems connected with semisimple Lie
  algebras. \emph{Invent. Math.} 37 (1976), 93--108.
\bibitem[OP2]{OP2} Olshanetsky, M.A., Perelomov, A.M., Quantum
  systems related to root systems, and radial parts of Laplace
  operators. \emph{Funct. Anal. Appl.} 12 (1978), 121--128. 
\bibitem[O1]{Op1} Opdam, E.M., Dunkl operators, Bessel functions and the 
discriminant of a finite Coxeter group. \emph{Compositio Math.} 85
(1993), 333--373.
\bibitem[O2]{Op} Opdam, E.M., Harmonic analysis for certain representations 
of graded Hecke algebras. \emph{Acta Math.} 175 (1995), 75 -- 121
\bibitem[Pa]{Pa} Pasquier, V.: A lecture on the Calogero-Sutherland models. 
In: Integrable models and strings (Espoo, 1993), \emph{Lecture Notes in Phys.}
 436,  Springer, 1994, pp. 36--48.
\bibitem[Pe]{Pe} Perelomov, A.M., Algebraical approach to the solution
  of  a one-dimensional model of $N$ interacting
  particles. \emph{Teor. Mat. Fiz.} 6 (1971), 364--391.
\bibitem[Po]{Po} Polychronakos, A.P., Exchange operator formalism for
  integrable systems of particles. \emph{Phys. Rev. Lett.} 69 (1992), 
703--705.
\bibitem[R1]{R1} R\"osler, M.,  Bessel-type signed hypergroups on $\b R$. 
In: Heyer, H., Mukherjea, A. (eds.) Probability measures on groups and 
related structures XI. Proceedings,  Oberwolfach 1994. 
World Scientific 1995,  292--304. 
\bibitem[R2]{R2} R\"osler, M.,  Generalized Hermite polynomials and the heat 
equation for Dunkl operators. \emph{Comm. Math. Phys.} 192 (1998),
519--542.
\bibitem[R3]{R3} R\"osler, M., Positivity of Dunkl's intertwining
  operator. \emph{Duke Math. J.} 98 (1999), 445--463.
\bibitem[R4]{R4} R\"osler, M., Short-time estimates for heat kernels associated with
root systems, in \textit{Special Functions,} Conf. Proc. Hong Kong June 1999, eds.
C. Dunkl et al., World Scientific, Singapore, 2000, 309--323.
\bibitem[R5]{R5} R\"osler, M., One-parameter semigroups related to abstract 
quantum models of Calogero type. In: Infinite Dimensional
              Harmonic Analysis (Kyoto, Sept. 1999, eds. H. Heyer et al.) 
Gr\"abner-Verlag 2000, 290--305. 
\bibitem[RdJ]{dJR} R\"osler, M., de Jeu, M., Asymptotic analysis for the Dunkl kernel.
 math.CA/0202083; to appear in J. Approx. Theory. 
\bibitem[RV]{RV3} R\"osler, M., Voit, M., Markov Processes related
  with Dunkl operators. \emph{Adv. Appl. Math.} 21 (1998), 575--643.
  \bibitem[Ros]{Ros} Rosenblum, M., Generalized Hermite polynomials and the 
Bose-like oscillator calculus. In: Operator Theory: Advances and 
Applications, Vol.  73,  Basel, Birkh\"auser Verlag 1994, 369--396.
\bibitem[Su]{Su} Sutherland, B., Exact results for a quantum many-body
  problem in one dimension. \emph{Phys. Rep.} A5 (1972), 1372--1376.
\bibitem[Ti]{Ti} Titchmarsh, E.C., \emph{The Theory of Functions.} 2nd ed.,
Oxford University Press, 1950.
\bibitem[T]{T} Trim\`eche, K., Paley-Wiener Theorems for the Dunkl transform
and Dunkl translation operators. \emph{Integral Transform. Spec. Funct.} 13 
(2002), 17--38.
\bibitem[W]{W} A. Wintner, On a theorem of B\^{o}cher in the theory of ordinary
 linear differential equations, \emph{Amer. J. Math.} 76 (1954), 183--190.
\bibitem[UW]{UW} Ujino, H., Wadati., M., Rodrigues formula for Hi-Jack
  symmetric polynomials associated with the quantum Calogero
  model. \emph{J. Phys. Soc. Japan} 65 (1996), 2423--2439.

\end{thebibliography}
\end{document}